\documentclass[11pt]{amsart}
\usepackage{amsmath,amsthm,amsfonts,amscd,amssymb}
\input psfig.sty
\newtheorem{theorem}{Theorem}[section]

\newtheorem{lemma}[theorem]{Lemma}
\newtheorem{corollary}[theorem]{Corollary}
\theoremstyle{definition}

\theoremstyle{remark}
\newtheorem{remark}[theorem]{Remark}
\numberwithin{equation}{section}

\DeclareMathOperator{\id}{id} 
 
 \DeclareMathOperator{\Area}{Area}
 \DeclareMathOperator{\Diff}{Diff}
\DeclareMathOperator{\Mob}{M\ddot{o}b}
\DeclareMathOperator{\Poin}{Poincar\acute{e}}
\DeclareMathOperator{\Mobius}{M\ddot{o}bius}
 \DeclareMathOperator{\PSL}{PSL}
\DeclareMathOperator{\PSU}{PSU}

\newcommand{\Z}{{\mathbb{Z}}}
\newcommand{\R}{{\mathbb{R}}}

\newcommand{\C}{{\mathbb{C}}}
\newcommand{\al}{\alpha}

\newcommand{\de}{\delta}
\newcommand{\g}{\gamma}
\newcommand{\s}{\gamma}
\newcommand{\si}{\sigma}

\newcommand{\bk}{\backslash}
\newcommand{\Ga}{\Gamma}

\newcommand{\w}{\overline{w}}

\newcommand{\pa}{\partial}
\newcommand{\la}{\langle}
\newcommand{\ra}{\rangle}
\newcommand{\qC}{\mathcal{C}}

\newcommand{\Ka}{K\"{a}hler\:}
\newcommand{\Te}{Teichm\"{u}ller\:}

\newcommand{\ov}{\overline}
\newcommand{\ep}{\epsilon}
\newcommand{\vep}{\varepsilon}

\newcommand{\z}{\bar{z}}

\newcommand{\ma}[4]{(\begin{smallmatrix}
              #1 & #2 \\ #3 & #4
             \end{smallmatrix})}

\begin{document}
\title{Velling-Kirillov metric on the universal Teichm\"uller curve}
\author{Lee-Peng Teo}
\email{lpteo@math.nctu.edu.tw}
\address{Department of Applied Mathematics, National Chiao Tung University, 1001,
 Ta-Hsueh Road, Hsinchu City, 30050, Taiwan R.O.C. }
\begin{abstract}
We extend Velling's approach and prove that the second variation of the spherical
areas of a family of domains
 defines a Hermitian metric on the universal Teichm\"{u}ller curve, whose pull--back
  to $\Diff+(S^1)/S^1$ coincides
 with the Kirillov metric. We call this Hermitian metric the Velling-Kirillov metric.
 We show that the vertical
 integration of the
 square of the symplectic form of the Velling-Kirillov metric on the universal
 Teichm\"uller curve is
 the symplectic form that
  defines the Weil-Petersson metric on the universal Teichm\"{u}ller space.
   Restricted to a finite dimensional
   Teichm\"uller space, the vertical integration of the corresponding form
   on the Teichm\"uller curve is also the
    symplectic form that defines the Weil-Petersson metric on the Teichm\"uller space.
\end{abstract}
\date{\today}
\maketitle

\section{introduction}
Let $T(1)$ be the universal \Te space and $\mathcal{T}(1)$ be the corresponding universal
\Te curve. $T(1)$ and $\mathcal{T}(1)$ have natural structure of infinite dimensional
 complex manifolds and the natural projection $p: \mathcal{T}(1) \rightarrow T(1)$
  ia a holomorphic fibration. In \cite{Ve}, J. Velling introduced a metric on
   $T(1)$ by using spherical areas. Namely, consider the Bers embedding of
    $T(1)$ into the Banach space
\[
 A_{\infty}(\Delta)=\left\{\phi \;\; \text{holomorphic on } \Delta:
  \sup_{z\in \Delta} |\phi(z)(1-|z|^2)^2| < \infty \right\},
\]
where $\Delta$ is the unit disc. For every $Q \in A_{\infty}(\Delta)$
and $t$ small, solution to the equation
\begin{align}\label{Sch}
\mathcal{S}(f^{tQ})= tQ,
\end{align}
where $\mathcal{S}(f)$ is the Schwarzian derivative of the function $f$,
defines a family of domains $\Omega_t = f^{tQ}(\Delta)$. Here $f^{tQ}$ is
normalized such that $f^{tQ}(0) = 0$, $f_z^{tQ}(0)=1$ and $f_{zz}^{tQ}(0)=0$.
J. Velling proved that the spherical area $A_S(\Omega_t)$ of the domain $\Omega_t$ satisfies
\[
\frac{d^2}{dt^2} A_S(\Omega_t)\bigr\vert_{t=0} \geq 0.
\]
This defines a Hermitian metric on the tangent space to $T(1)$ at the origin,
 identified with $A_{\infty}(\Delta)$. Our first result,
 Theorem \ref{sphericalmetric}, is the following explicit formula for this metric:
\begin{align*}
\parallel Q \parallel_S^2 &= \frac{1}{2 \pi} \frac{d^2}{dt^2}
A_S(\Omega_t)\bigr\vert_{t=0} = \sum_{n=2}^{\infty} n |a_n|^2,
\end{align*}
where $Q(z) = \sum_{n=2}^{\infty} (n^3-n) a_n z^{n-2}$.
The series converges for all $Q \in A_{\infty}(\Delta)$.
 However, since the spherical area of the domain $f^{tQ}(\Delta)$ is not
  independent of the choice of the function $f^{tQ}$ that satisfies
  \eqref{Sch}, $\parallel \cdot \parallel_S$ does not naturally define
   a metric on $T(1)$ by right group translations\footnote{The metric on $T(1)$
    defined as a pull-back of the Hermitian metric on $A_{\infty}(\Delta)$
     given by $\parallel \cdot \parallel_S$ is not natural.
      It does not induce a metric on finite dimensional Teichm\"uller spaces
       embeded in $T(1)$ since these embeddings are base point dependent.}.
        Nevertheless, we observe that Velling's approach can be generalized
         to define a metric on the universal \Te curve $\mathcal{T}(1)$.
          This is achieved by a natural identification of $\mathcal{T}(1)$
           with the space $\text{Homeo}_{qs}(S^1)/S^1$ --- the subgroup of
            orientation preserving quasisymmetric homeomorphisms of
             the unit circle that fix the point $1$, and with the space
\begin{align*}
\tilde{\mathcal{D}} =\{ f:\Delta \longrightarrow \hat{\C}\; \text{a univalent function}
 \; : f(0)=0, f'(0)=1, \\ f \; \text{has a
quasiconformal extension to}\; \hat{\C} \},
\end{align*}
which we prove in Section 2. This endows $\mathcal{T}(1)$ with a
group structure \footnote{It is well known (see, e.g., \cite{Nag,
Lehto}) that $\mathcal{T}(1)$ is not a topological group.}.
Following Velling's approach to $T(1)$, given a one-parameter
family of univalent functions $f^t : \Delta \rightarrow \hat{\C}
\in \tilde{\mathcal{D}}$, $f^t \vert_{t=0} =\id$, which defines a
tangent vector $\upsilon$ corresponding to $\frac{d}{dt} f^t
\vert_{t=0}$ at the origin, we define a metric on the tangent
space to $\mathcal{T}(1)$ at the origin by
\[
\parallel \upsilon \parallel^2 = \frac{1}{2\pi}\frac{d^2}{dt^2}
A_S(f^{t}(\Delta))\bigr\vert_{t=0},
\]
and extend it to every point of $\mathcal{T}(1)$ by right
translations. This metric is Hermitian and \Ka. More remarkably,
its pull-back via the embedding $\Diff_+(S^1)/S^1 \hookrightarrow
\text{Homeo}_{qs}(S^1)/S^1 \simeq \mathcal{T}(1)$ is precisely the
 metric
\[
\parallel \upsilon \parallel^2 = \sum_{n=1}^{\infty} n |c_n|^2
\]
on $\Diff_+(S^1)/S^1$ introduced by Kirillov \cite{Ki, KY} via the
coadjoint orbit method. Here $\upsilon = \sum_{n} c_n e^{i n
\theta} \frac{\pa}{\pa \theta}$, $c_{-n} = \ov{c_n}$ is a vector
field on $S^1$. We call this K\"ahler metric on $\mathcal{T}(1)$
the Velling-Kirillov metric and prove that it is the unique right
invariant K\"ahler metric on $\mathcal{T}(1)$.

Let $\kappa$ be the symplectic form of the Velling-Kirillov metric
on $\mathcal{T}(1)$. We consider the $(1,1)$ form $\omega$ on
$T(1)$, which is the vertical integration of the $(2,2)$ form
$\kappa \wedge \kappa$ on $\mathcal{T}(1)$, i.e., integration of
$\kappa \wedge \kappa$ over the fibers of the fibration $p :
\mathcal{T}(1) \rightarrow T(1)$. We proved that this is
equivalent to Velling's suggestion of averaging the Hermitian form
$\parallel \cdot \parallel_S$ along the fibers of $\mathcal{T}(1)$
over $T(1)$. Our second result, which we prove in Theorem
\ref{averagecoefficient} and Theorem \ref{Vellingmetric}, is that
$\omega$ is the symplectic form of the Weil-Petersson metric on
$T(1)$, defined only on tangent vectors which correspond to
$H^{\frac{3}{2}}$ vector fields on $S^1$.

When $\Gamma$ is a cofinite Fuchsian group, the \Te space of
$\Gamma$, $T(\Gamma)$ embeds holomorphically in
$T(1)$. The Bers fiber space $\mathcal{BF}(\Gamma)$ is the
inverse image of $T(\Gamma)$ under the
 projection map $\mathcal{T}(1) \rightarrow T(1)$ and the
 Teichm\"uller curve $\mathcal{F}(\Gamma)$ is
  a quotient space of $\mathcal{BF}(\Gamma)$. The symplectic
   form $\kappa$ is well defined when restricted
  to $\mathcal{F}(\Gamma)$. We prove in Theorems \ref{average2}
   and \ref{Velling2} that the vertical integration of $\kappa \wedge \kappa$
   via the map $\mathcal{F}(\Gamma) \rightarrow T(\Gamma)$ is the symplectic
   form that defines the Weil-Petersson metric on $\Gamma$.

In the appendix, we consider an analog of Bers embedding for $\mathcal{T}(1)$. We
prove that $\mathcal{T}(1)$ embeds into the Banach space
\begin{align*}
\mathcal{A}_{\infty}(\Delta) =\left\{\psi \;\text{holomorphic on }
\Delta: \sup_{z\in \Delta} |\psi(z)(1-|z|^2)| < \infty \right\},
\end{align*}
and its image contains an open ball about the origin of
$\mathcal{A}_{\infty}(\Delta)$. We also verify that $\mathcal{A}_{\infty}(\Delta)$
 and $A_{\infty}(\Delta) \oplus \C$ induce the same complex structure on
  $\mathcal{T}(1)$. These results are not used in the main text.

The content of this paper is the following. In Section 2, we review different
 models for
 the universal \Te space, the universal \Te curve and study
 their relations with the homogenuous spaces of
 $\text{Homeo}_{qs}(S^1)$. In Section 3, we review
 Velling's approach and define a metric on the universal \Te curve.
  We prove that its pull-back to $\Diff_+(S^1)/S^1$ coincides
  with the Kirillov metric. In Section 4, we prove that
   the vertical integration of the square of the symplectic
   form of Velling-Kirillov metric is the symplectic form
   that defines the Weil-Petersson metric on \Te spaces. In
   the appendix, we consider an embedding of $\mathcal{T}(1)$.

\vspace{0.2cm}
\noindent
\textbf{Acknowledgements.} This work is an extension of a part of my
Ph.D. thesis. I am especially grateful to
 my advisor Leon A. Takhtajan for the stimulating discussions and useful suggestions.
  I would also like to thank him for bringing this subject to my
  attention. J.Velling has kindly made his unpublished manuscript
  \cite{Ve} available, which has been a great stimulation for
  the present work. The author has quoted or reproduced some of his results for
  the convenience of the reader.

\section{Universal \Te Space and Universal \Te Curve}
\subsection{\Te theory}
Here we collect basic facts from \Te theory. For details, see \cite{Nag2, A2, Lehto}.

Let $T(1)$ be the universal \Te space. There are two classical models of this space.

Let $\Delta$ be the open unit disc, $\Delta^*$ be the exterior of the unit disc.
Let $L^{\infty}(\Delta^*)$ (resp. $L^{\infty}(\Delta)$) be the complex Banach
space of bounded Beltrami differentials on $\Delta^*$ (resp. $\Delta$) and let
$L^{\infty}(\Delta^*)_1$ be the unit ball of $L^{\infty}(\Delta^*)$.
For any $\mu \in L^{\infty}(\Delta^*)_1$, we consider the following two constructions.
\begin{enumerate}
\item
Model A: $w_{\mu}$ theory.

\noindent
We extend $\mu$ by reflection to $\Delta$, i.e.
\begin{equation}\label{sym}
\mu(z) = \ov{\mu\left(\frac{1}{\z}\right)} \frac{z^2}{\z^2}\; , \hspace{2cm} z \in \Delta.
\end{equation}
There is a unique quasiconformal map $w_{\mu}$, fixing $-1, -i$ and $1$,
which solves the Beltrami equation
\begin{align*}
  (w_{\mu})_{\z} &= \mu (w_{\mu})_z \;.
\end{align*}
It satisfies
\begin{align}\label{reflection}
  w_{\mu} (z) &= \ov{w_{\mu} \left(\frac{1}{\z}\right)}
\end{align}
due to the reflection symmetry ~\eqref{sym}. As a result,
$w_{\mu}$ fixes the unit circle $S^1$, $\Delta$ and $\Delta^*$.

\vspace{0.2cm}
\item
Model B: $w^{\mu}$ theory.

\noindent
We extend $\mu$ to be zero outside $\Delta^*$. There is a unique
 quasiconformal map $w^{\mu}$, holomorphic on the unit disc, which
  solves the Beltrami equation
\begin{equation*}
   w_{\z}^{\mu} = \mu w_z^{\mu},
\end{equation*}
and is normalized such that $f = w^{\mu} \vert_{\Delta}$ satisfies
$f(0)=0$, $f'(0)=1$ and $f''(0)=0$.
\end{enumerate}

The universal \Te space $T(1)$ is defined as a set of equivalence
 classes of normalized quasiconformal maps
\[
  T(1) = L^{\infty}(\Delta^*)_1/\sim ,
\]
where $\mu\sim \nu$ if and only if $w_{\mu} = w_{\nu}$ on the unit
circle, or equivalently, $w^{\mu}= w^{\nu}$ on the unit disc.

Using model B, we can identify $T(1)$ with the space
\begin{align*}
\mathcal{D} = \{ f: \Delta \rightarrow \C \; \text{univalent} \; :
f(0)=0, f'(0)=1, f''(0) = 0;\\
 f \;\text{has a quasiconformal extension to}\; \C\}.
\end{align*}

Let $\mathcal{S}(f)$ be the Schwarzian derivative of the
function $f$, which is given by
\[
\mathcal{S}(f) = \left(\frac{f_{zz}}{f_z}\right)_z - \frac{1}{2}
 \left(\frac{f_{zz}}{f_z}\right)^2 = \frac{f_{zzz}}{f_z} -
 \frac{3}{2}\left(\frac{f_{zz}}{f_z}\right)^2 .
\]
Let $A_\infty(\Delta)$ be the Banach space
\[
A_\infty(\Delta)= \left\{\phi \;\text{holomorphic on }
\Delta: \sup_{z\in \Delta} |\phi(z)(1-|z|^2)^2| < \infty \right\}.
 \]
The Bers embedding $T(1) \hookrightarrow A_{\infty}(\Delta)$ which
 maps $[\mu]$ --- the equivalence class of $\mu$ --- to
 $\mathcal{S}(w^{\mu}|_{\Delta})$
endows $T(1)$ with a unique structure of a complex Banach
 manifold such that the projection map
\[
\Phi : L^{\infty}(\Delta^*)_1 \rightarrow T(1)
\]
is a holomorphic submersion. In particular, $L^{\infty}(\Delta^*)_1$
and $A_{\infty}(\Delta)$ induce the same complex structure on $T(1)$.

The derivative of the map $\Phi$ at the origin
\[
D_0 \Phi : L^{\infty}(\Delta^*) \longrightarrow T_0 T(1)
\]
is a complex linear surjection, with kernel $\mathcal{N}(\Delta^*)$
--- the space of infinitesimally trivial Beltrami differentials. Explicitly,
\begin{align*}
\mathcal{N}(\Delta^*) = \left\{ \mu \in L^{\infty}(\Delta^*) :
\iint\limits_{\Delta^*} \mu \phi = 0\;, \; \forall \phi \in
A_1(\Delta^*)\right\} \;,
\end{align*}
where $A_1(\Delta^*)$ is the Banach space of $L^1$
(with respect to Lebesgue measure on $\Delta^*$)
 holomorphic functions on $\Delta^*$.

Define
\[
A_\infty(\Delta^*)= \left\{\phi \;\text{holomorphic on }
\Delta^*: \sup_{z\in \Delta^*} |\phi(z)(1-|z|^2)^2| < \infty \right\}
\]
and its complex anti-linear isomorphic space
\[
\Omega^{-1,1}(\Delta^*)= \left\{ \mu(z) = (1-|z|^2)^2 \ov{\phi(z)}:
\phi \in A_{\infty} (\Delta^*)\right\},
\]
the space of harmonic Beltrami differentials on $\Delta^*$.
There is a canonical splitting
\[
L^{\infty}(\Delta^*) = \mathcal{N}(\Delta^*) \oplus \Omega^{-1,1}(\Delta^*),
\]
which identifies the tangent space at the origin of $T(1)$
with $\Omega^{-1,1}(\Delta^*)$. Moreover, Bers embedding
 induces the isomorphism $\Omega^{-1,1}(\Delta^*)
 \xrightarrow{\sim} A_{\infty}(\Delta)$ given by
\begin{align}\label{Bersmap}
\mu \mapsto \phi(z) = -\frac{6}{\pi} \iint\limits_{\Delta^*}
\frac{\mu(\zeta)}{(\zeta-z)^4} \left\vert \frac{d\zeta \wedge
d\ov{\zeta}}{2}\right\vert.
\end{align}

$L^{\infty}(\Delta^*)_1$ has a group structure induced by
the composition of quasiconformal maps,
\begin{align*}
\lambda * \mu = \nu\;, \hspace{1cm} \text{where} \hspace{0.5cm}
w_{\nu} = w_{\lambda}\circ w_{\mu}\; .
\end{align*}
Explicitly, it is given by
\begin{align*}
\nu = \frac{ \mu + (\lambda\circ w_{\mu})\frac{\ov{(w_{\mu})_z}}
{(w_{\mu})_z}}{1+ \ov{\mu}(\lambda\circ w_{\mu})
\frac{\ov{(w_{\mu})_z}}{(w_{\mu})_z}}\; .
\end{align*}
This group structure descends to $T(1)$. Moreover, the right
group translation by $[\mu]$, $R_{[\mu]} : T(1) \rightarrow T(1)$,
$[\lambda] \mapsto [\lambda* \mu]$ is biholomorphic.
However, the left group translation is not even a continuous map
on $T(1)$ (see, e.g. \cite{Nag2, Lehto}).

\begin{remark}
Conventionally, the model of the universal \Te space is the
complex conjugate of the one we define above. Consider the natural
complex anti-linear isomorphism
\begin{align*}
 L^{\infty}(\Delta^*)_1 &\hspace{0.5cm} \rightarrow
 \hspace{0.5cm} L^{\infty}(\Delta)_1 \\
 \mu &\hspace{0.5cm} \mapsto \hspace{0.5cm} \tilde{\mu} =
  \ov{\mu\left(\frac{1}{\z}\right)} \frac{z^2}{\z^2}\; ,
  \hspace{0.5cm} z \in \Delta \;.
\end{align*}
Setting $\tilde{\mu}$ to be zero outside $\Delta$,
we obtain a unique solution of the Beltrami equation
\[
w_{\z}^{\tilde{\mu}} = \tilde{\mu} w_z^{\tilde{\mu}},
\]
which is holomorphic on $\Delta^*$ and normalized such that
$g = w^{\tilde{\mu}}\vert_{\Delta^*}$ has Laurent expansion at $\infty$ given by
\begin{align}\label{Laurent}
g(z) = z ( 1+ \frac{a_2}{z^2} + \frac{a_3}{z^3} + \ldots).
\end{align}
Thus $T(1)$ is identified with the space
\begin{align*}
\mathcal{D^*} = \{ g: \Delta^* \rightarrow \hat{\C} \;
\text{univalent} \; : g \; \text{has Laurent expansion at
$\infty$ given by \eqref{Laurent}}  \\
 \text{and has quasiconformal extension to}\; \hat{\C\}}.
\end{align*}
\end{remark}

The universal \Te curve $\mathcal{T}(1)$ is a fiber space over $T(1)$.
The fiber over each point $[\mu]$ is the quasidisc $w^{\mu}(\Delta^*)
\in \hat{\C}$ with the complex structure induced from $\hat{\C}$,
\begin{align}\label{universalcurve}
\mathcal{T}(1) = \Bigl\{ ([\mu], z) : \; [\mu] \in T(1), \;z
\in w^{\mu}(\Delta^*)\Bigr\}.
\end{align}
It is a Banach manifold modeled on $A_{\infty}(\Delta) \oplus \C$
\footnote{Here we give $\hat{\C} \setminus \{0\}$ the complex Banach
manifold structure of $\C$ via the conformal map
$z \mapsto \frac{1}{z}$.}. We have a real analytic isomorphism
between $ T(1) \times \Delta^* $ and $\mathcal{T}(1)$ given by
\[
([\mu], z ) \mapsto ([\mu], w^{\mu}(z)).
\]

\subsubsection{Homogenuous spaces of $\text{Homeo}_{qs}(S^1)$}
Let $\text{Homeo}_{qs}(S^1)$ be the group of orientation
preserving quasisymmetric homeomorphisms of the unit circle $S^1$.
It contains the subgroup of orientation preserving diffeomorphisms
--- $Diff_+(S^1)$. We denote by $\Mob(S^1)$ the subgroup of
$\Mobius$ transformations and abusing notations, denote by $S^1$
the subgroup of rotations.

Consider the model A of the universal \Te space $T(1)$ given
above. Clearly, the map $T(1) \ni [\mu] \mapsto
w_{\mu}\vert_{S^1}\in \text{Homeo}_{qs}(S^1)$ is well defined and
is one-to-one. Ahlfors-Beurling extension theorem implies that its
image consists of all normalized orientation preserving
quasisymmetric homeomorphisms of the unit circle (see, e.g.,
\cite{LBers, Nag2, Lehto}), in other words,
\[
T(1) \cong \text{Homeo}_{qs}(S^1)/\Mob(S^1).
\]

Let $\mu\in \Omega^{-1,1}(\Delta^*)$ be a tangent vector
at the origin of $T(1)$. It generates the one-parameter flow $w_{t\mu}$
 and the corresponding vector field is given by $\dot{w}_{\mu}
 \frac{\pa}{\pa z}$, where
\begin{align*}
\dot{w}_{\mu}(z) &= \frac{(z+1)(z+i)(z-1)}{2\pi i} \iint\limits_{\C}
\frac{\hat{\mu}(\zeta)}{(\zeta - z)(\zeta+1)(\zeta+i)(\zeta-1)}
d\zeta \wedge d \ov{\zeta},
\end{align*}
and $\hat{\mu}$ is the extension of $\mu$ by reflection to $\C$.
Restricted to $S^1$, we have $\dot{w}_{\mu}(z) = iz\mathfrak{u}(z)$,
where $\mathfrak{u}(e^{i\theta}) \frac{\pa}{\pa \theta}$ is the
vector field on $S^1$.

It was proved by Reimann (see \cite{Reimann,GaSu,Nag}) that the
tangent space to $\text{Homeo}_{qs}(S^1)$ at the origin is the Zygmund space
\begin{align*}
\Lambda(S^1)= \Biggl\{\mathfrak{u}(e^{i\theta}) \frac{\pa}{\pa \theta}
: (i)&\mathfrak{u} : S^1 \rightarrow \R \; \text{is continuous} . \\
 (ii) & F_{\mathfrak{u}}(x) =\frac{1}{2}(x^2+1) \mathfrak{u}
 \left(\frac{x-i}{x+i}\right) \; \text{is in} \;\; \Lambda(\R).\Biggr\} ,
\end{align*}
where
\begin{align*}
\Lambda(\R) = \Biggl\{ &F: \R \rightarrow \R: (i) F \; \text{is continuous}.\\
&(ii) \left\vert F(x+t)  + F(x-t) - 2F(x) \right\vert \leq B
\left\vert t \right\vert \; \text{for some} \;\; B,  \forall x, t \in \R \Biggr\}.
\end{align*}
By imposing extra normalization conditions, we can characterize the
tangent space at the origin of $\text{Homeo}_{qs}(S^1)/S^1$ and
$\text{Homeo}_{qs}(S^1)/\Mob(S^1)$ in a similar way.

\begin{remark}
It is not known how to characterize the Zygmund space
$\Lambda(S^1)$ using Fourier coefficients on $S^1$.
\end{remark}

In \cite{Ki}, Kirillov considered the Lie group $\Diff_+(S^1)$ and
proved that there is a natural bijection between the space $\mathcal{K}$
of smooth contours of conformal radius $1$ which contain $0$
in their interior and the space $\Diff_+(S^1)/S^1$.
We generalize this bijection in the following theorem.

\begin{theorem}
There is a natural bijection between the space $\text{Homeo}_{qs}(S^1)/S^1$
and the space $\mathcal{K}_{qc}$ of all quasicircles, i.e. images of the unit
circle under a quasiconformal maps, of conformal radius 1 which contain $0$
in their interior. Moreover, for every $\g \in \text{Homeo}_{qs}(S^1)/S^1$,
there exists two univalent functions $f: \Delta \rightarrow \C$ and
$g: \Delta^* \rightarrow \hat{\C}$ determined by the following properties:
\begin{enumerate}
\item
$f$ and $g$ admit quasiconformal extensions to quasiconformal
mappings of $\hat{\C}$.
\item
$ \g= g^{-1} \circ f\vert_{S^1} \mod S^1.$
\item
$f(0)=0, \; f'(0)=1$.
\item
$g(\infty) = \infty, \; g'(\infty) > 0$.

\end{enumerate}
\end{theorem}
\begin{proof}
By Ahlfors-Beurling extension theorem, an orientation preserving
quasisymmetric homeomorphism $\g$ of the unit circle can be extended to
a quasiconformal map $w$ of $\hat{\C}$ satisfying the reflection property
\eqref{reflection}.
Let $\mu$ be the Beltrami differential of the map $w\vert_{\Delta^*}$.
Up to a linear fractional transformation, $w$ agrees with $w_{\mu}$
that we define in Section 2.1, i.e. $w = \si_1 \circ w_{\mu}$ for
some $\si_1 \in \PSU(1,1)$. The corresponding map $w^{\mu}$ (Section 2.1)
is holomorphic inside the unit disc $\Delta$. Define
$g = \si_2 \circ w^{\mu} \circ w^{-1}$,
where $\si_2 \in \PSL(2, \C)$ is uniquely determined by the
conditions $f = \si_2 \circ w^{\mu}$ satisfies $f(0)=0$,
$f'(0)=1$ and $g$ satisfies $g(\infty)= \infty$.
The maps $f\vert_{\Delta}$ and $g\vert_{\Delta^*}$ are holomorphic.
They do not depend on the extension of $\g$ and we have $\g = g^{-1}
\circ f \vert_{S^1}$. The image of $S^1$ under $f$, which is the same
as the image of $S^1$ under $g$, is by definition a quasicircle $\qC$
with conformal radius $1$. By post-composing $w$ with a rotation,
the map $g$ also satisfies $g'(\infty) >0$.

Conversely, by definition a quasicircle $\qC$ is the image of $S^1$
under a quasiconformal map $h: \C \rightarrow \C$. Let $\mu_1$ be the
Beltrami differential of $h\vert_{\Delta}$, extended to $\Delta^*$
by reflection. Let $w_{\mu_1}$ be a solution of the corresponding
Beltrami equation. Then $f=h \circ w_{\mu_1}^{-1}$ is a
quasiconformal map that is holomorphic inside $\Delta$.
When $0$ is in the interior of $\qC$, there is a unique way to
normalize $w_{\mu_1}$ by post-composition with a $\PSU(1,1)$
transformation such that $f(0)=0$ and $f'(0)>0$. The image of $S^1$
under $f$ is the quasicircle $\qC$. In fact by Riemann mapping theorem,
$f\vert_{\Delta}$ is uniquely determined by $\qC$ and the
normalization conditions $f(0)=0$, $f'(0)>0$. $\qC$ has
conformal radius $1$ implies that $f'(0)=1$. Let $\mu$ be the
Beltrami differential of $f\vert_{\Delta^*}$, extended to $\Delta$
by reflection. Let $w_{\mu}$ be a solution of the corresponding
Beltrami equation. Define $g = f \circ w_{\mu}^{-1} \circ \si$,
where $\si \in \PSU(1,1)$ is uniquely determined so that
$g(\infty) = \infty$ and $g'(\infty) >0$. The map $\g = g^{-1}
\circ f\vert_{S^1}$ is then an orientation preserving
quasisymmetric homeomorphism of the unit circle.
\end{proof}
The decomposition $\g = g^{-1}\circ f$ is known as
conformal welding.
 Using the fact that the correspondence between $f$ and the
 quasicircle $\qC$ is one-to-one, we can identify
 $\text{Homeo}_{qs}(S^1)/S^1$ with the space of univalent functions
\begin{align*}
\tilde{\mathcal{D}} = \{ f:\Delta \longrightarrow \hat{\C}\;\,
\text{a univalent function} \; : f(0)=0, f'(0)=1, \\ f \; \text{has a
quasiconformal extension to}\; \C \}.
\end{align*}
$\tilde{\mathcal{D}}$ is a complex subspace of the complex space
of sequences $\{a_n\}$ (Fourier coefficients of the holomorphic
function $f$). This induces a complex structure on
$\text{Homeo}_{qs}(S^1)/S^1$.

\begin{remark}\label{equivalenta}
Notice that if $\g = w_{\mu}\vert_{S^1}$ up to post-composition with a
$\PSU(1,1)$ transformation, then corresponding $f$ is equal to
$w^{\mu}$ up to post-composition with a $\PSL(2, \C)$ transformation.
\end{remark}


We identify $\text{Homeo}_{qs}(S^1)/S^1$ as a subgroup of
$\text{Homeo}_{qs}(S^1)$ consisting of quasisymmetric homeomorphisms
that fix the point $1$. Consider $\text{Homeo}_{qs}(S^1)/$ $ \Mob(S^1)$
as the subspace of $\text{Homeo}_{qs}(S^1)/S^1$ corresponding to
the natural inclusion $T(1) \simeq \mathcal{D} \hookrightarrow
\tilde{\mathcal{D}} \simeq \text{Homeo}_{qs}(S^1)/S^1$.
Analogous to the isomorphism $T(1) \simeq
\text{Homeo}_{qs}(S^1)/\Mob(S^1)$, we have
\begin{theorem}\label{universal}
There is an isomorphism between $\mathcal{T}(1)$ and
$\text{Homeo}_{qs}(S^1)/S^1 \simeq \tilde{\mathcal{D}}$.
Moreover, the complex structure of $\mathcal{T}(1)$
induced from $A_{\infty}(\Delta) \oplus \C$ coincides
with the complex structure induced from $\tilde{\mathcal{D}}$.
\end{theorem}
\begin{proof}
 The fiber of $\text{Homeo}_{qs}(S^1)/S^1$ over $\g \in
 \text{Homeo}_{qs}(S^1)/\Mob(S^1)$, consists of all quasisymmetric
 homeomorphisms of the form $\si \circ \g$ mod $S^1$, where
 $\si \in \PSU(1,1)$ mod $S^1$ are parametrized by
 $w \in \Delta^* \simeq \PSU(1,1)/S^1$, i.e.
\begin{align}\label{sigmamap}
\si_w(z) = \frac{1 - z \bar{w}}{z-w}.
\end{align}
Let $f, g$ (resp. $f_w, g_w$) be the univalent functions
corresponding to $\g$ (resp. $\g_w = \si_w \circ \g$), i.e.
\begin{align*}
\g = g^{-1}\circ f, \hspace{2cm} \si_w \circ \g = g_w^{-1} \circ f_w.
\end{align*}
Using Remark \ref{equivalenta}, we have
\begin{align*}
f_w = \lambda_w \circ f, \hspace{1cm}\text{and hence}
\hspace{0.5cm}  g_w = \lambda_w \circ g \circ \si_w^{-1}
\end{align*}
for some $\lambda_w \in \PSL(2, \C)$. The normalization conditions
on $f_w \in \tilde{\mathcal{D}}$ and $f \in \mathcal{D}$ imply that
\begin{align}\label{lambda}
\lambda_w (z)= \frac{z}{c_w z+1},  \hspace{0.5cm} \text{where}
\;\; c_w = -\frac{1}{2} \frac{f_w''(0)}{f_w'(0)}.
\end{align}
The condition $g_w(\infty) = \infty$ implies that
\begin{align}\label{constantc}
c_w = - \frac{1}{g(w)}.
\end{align}
Let $[\mu]$ be the equivalence class which corresponds to $\g$
under the isomorphism $T(1) \simeq \text{Homeo}_{qs}(S^1)/
\Mob(S^1)$. For $w \in \Delta^*$, since $f(\Delta^*) =
g(\Delta^*)$, the point $g(w)$ lies in $f(\Delta^*)=
w^{\mu}(\Delta^*)$. Hence the natural correspondence between
$\text{Homeo}_{qs}(S^1)/S^1$ ($\simeq \tilde{\mathcal{D}}$) and
$\mathcal{T}(1)$ given by
\begin{align}\label{isomor}
\si_w \circ \g \;\; ( f_w = \lambda_w \circ f) \mapsto ([\mu], g(w)),
\end{align}
is an isomorphism.

In the identification above, $T(1)$ is the natural
subspace $\{([\mu], \infty) : [\mu] \in T(1)\}$ of $\mathcal{T}(1)$.
The embedding $([\mu], \infty) \mapsto f$ of $T(1)$ into
$\tilde{\mathcal{D}}$ is the pre-Bers embedding.
Hence the complex structure of $T(1) \simeq A_{\infty}(\Delta)$
agrees with the complex structure induced from $\tilde{\mathcal{D}}$ .
From \eqref{isomor}, \eqref{lambda}, \eqref{constantc}, we see that
if we fix $[\mu]$ in $([\mu], z)\in \mathcal{T}(1)$, and change $z$
holomorphically, the corresponding $f\in \tilde{\mathcal{D}}$ associated
to $([\mu], z)$ changes by post-composition with
$\lambda = \ma{1}{0}{c}{1}\in \PSL(2, \C)$, where the coefficient
$c$ depends holomorphically on $z$. This implies that
the complex structure of $\mathcal{T}(1)$ induced from the embedding
$\mathcal{T}(1) \hookrightarrow A_{\infty}(\Delta) \oplus \C$ agrees
with the complex structure of $\tilde{\mathcal{D}}$
induced from the isomorphism \eqref{isomor}.

\end{proof}
We can identify each point in $\mathcal{T}(1)$ as equivalence
classes of quasiconformal mappings as in the proof of the theorem above.
This immediately implies that $\mathcal{T}(1)$ also has a
group structure coming from composition of quasiconformal maps,
which is an extension of the group structure on $T(1)$.
Using the definition and the identification given in the proof
of Theorem \ref{universal}, the group multiplication in terms
of coordinates \eqref{universalcurve} is given by
\begin{align}
([\lambda], z)* ([\mu], z_0) &= ([\nu], z'), \label{gp1}
\end{align}
where
\begin{align}
\nu = \frac{ \mu + (\lambda\circ w)\frac{\ov{w_z}}{w_z}}{1+ \ov{\mu}
(\lambda\circ w)\frac{\ov{w_z}}{w_z}} \hspace{1cm} &\text{and}
\hspace{0.5cm} z' = w^{\nu} \circ w^{-1} \circ (w^{\lambda})^{-1}(z) \label{gp2}.
\end{align}
Here $w$ is the quasiconformal map corresponding to the point
$([\mu], z_0)\in \mathcal{T}(1) \simeq
\text{Homeo}_{qs}(S^1)/S^1$. The right group translation by
$([\mu], z_0)$, $R_{([\mu], z_0)}: \mathcal{T}(1) \rightarrow
\mathcal{T}(1)$ is biholomorphic (see \cite{Bersfiber}). Thus we
can identify the tangent space at $([\mu], z_0)$ with the tangent
space at $(0, \infty)$ --- the origin of $\mathcal{T}(1)$ --- via
the inverse of the derivative of the map $R_{([\mu], z_0)}$ at the
origin, i.e. via the map $\left( D_{(0,\infty)} R_{([\mu],
z_0)}\right)^{-1}$ . Moreover, this identification and the group
structure give rise to a splitting of the tangent space at each
point of $\mathcal{T}(1)$ into horizontal and vertical directions.
At the origin $(0, \infty)$, the vertical direction is spanned by
$\{0\} \oplus \C$ and the horizontal direction is spanned by
$\Omega^{-1,1}(\Delta^*)\oplus \{0\}$. A horizontal vector $(\nu,
0)$, $\nu \in \Omega^{-1,1}(\Delta^*)$ at the origin $(0, \infty)$
has a unique horizontal lift to each point $(0, z)$ on the fiber
at $(0, \infty)$. Namely, let $([t \nu], z_t')$, $z_0'=z $ be a
curve that defines the horizontal lift of $(\nu, 0)$ at the point
$(0, z)$. For $t$ small, $z_t'$ is determined by the following
equation
\begin{align*}
([\lambda(t)], \infty) * (0, z) = ([t \nu], z_t') , \hspace{1cm}
\lambda(t) \in  L^{\infty}(\Delta^*)_1.
\end{align*}
The point $(0,z)$ corresponds to the map $\si _z$ defined by
\eqref{sigmamap} (the subscript $z$ is not derivative).
Using the formulas \eqref{gp1}, \eqref{gp2}, taking derivative
with respect to $t$ and setting $t=0$ (which we denote by \;$\cdot$\;),
we have
\begin{align}\label{vectoridentification}
 \dot{\lambda} = \left( \nu \;\; \frac{\si_z'}{\ov{\si_z'}}\right)
 \circ \si_z^{-1} \hspace{1cm} \text{and} \hspace{1cm} \dot{ z}'
 = \dot{ w}^{ \nu} (z).
\end{align}
Hence the horizontal tangent vector $(\nu, 0)$ at $(0, \infty)$
is lifted to the vector $(\nu, \dot{ w}^{ \nu} (z))$ at $(0, z)$
and the latter is identified with the horizontal tangent vector
$ \left(\dot{\lambda} , 0 \right)$ at the origin $(0, \infty)$
of $\mathcal{T}(1)$.

\subsection{Identification of tangent spaces}
The isomorphism
\begin{align*}
\mathcal{W} : \text{Homeo}_{qs}(S^1)/S^1 &\longrightarrow
\tilde{\mathcal{D}},\\
\g & \mapsto f
\end{align*}
establishes the relation between real analytic
(through $\text{Homeo}_{qs}(S^1)/S^1)$) and complex analytic
(through $\mathcal{\tilde{D}}$) descriptions of $\mathcal{T}(1)$.
Infinitesimally, it takes the following explicit form.

\begin{theorem}\label{W}
The derivative of $\mathcal{W}$ at the origin is the linear mapping
$D_0 \mathcal{W} : T_0 \text{Homeo}_{qs}(S^1)/S^1 \rightarrow T_0
\tilde{\mathcal{D}}$ given by
\begin{align*}
\sum_{n \neq 0} c_n e^{in \theta} &\mapsto i \sum_{n=1}^{\infty} c_n z^{n+1}.
\end{align*}
\end{theorem}
\begin{proof}
Consider the smooth one parameter flow $\gamma^t = (g^t)^{-1} \circ f^t
\vert_{S^1}$, $\g^t \vert_{t=0} = \id$. It is known (see, e.g., \cite{Lehto})
that $\g^t$, $f^t$ and $g^t$ can be extended to quasiconformal mappings of
$\hat{\C}$, real analytic on $\hat{\C} \setminus S^1$.
The corresponding vector fields
\[
\frac{d}{dt} \g^t, \hspace{1.5cm} \frac{d}{dt}f^t, \hspace{0.8cm} \text{and}
\hspace{0.8cm} \frac{d}{dt} g^t
\]
are continuous on $\hat{\C}$, real analytic on $\hat{\C}\setminus S^1$.

We write the perturbative expansion
\begin{align*}
f^t(z) & = z + tu + O(t^2) = z + t z( a_1 z + a_2 z^2 + \ldots) + O(t^2),
\end{align*}
for $z \in \Delta$ and
\begin{align*}
g^t(z) &= z + tv + O(t^2) = z + t z (b_0 + b_1 z^{-1} + b_2 z^{-2} + \ldots)+O(t^2),
\end{align*}
for $z \in \Delta^*$.

We denote by
\[
\dot{\g} = \frac{d}{dt} \g^t \Bigr\vert_{t=0}, \hspace{1cm} \dot{f} =
\frac{d}{dt} f^t \Bigr\vert_{t=0} \hspace{0.5cm} \text{and}
\hspace{0.5cm} \dot{g} = \frac{d}{dt} g^t \Bigr\vert_{t=0},
\]
so that $\dot{f}\vert_\Delta = u$ and $\dot{g}\vert_{\Delta^*} =v$.

Under the Bers embedding, $\mathcal{S}(f^t\vert_\Delta)$ belongs to
a bounded subspace of $A_{\infty}(\Delta)$ and the corresponding
tangent vector to $T(1)$ at the origin is
\[
u_{zzz} = \frac{d}{dt}\mathcal{S}(f^t\vert_\Delta)\Bigr\vert_{t=0}
\in A_{\infty}(\Delta).
\]

To continue the proof, we need the following two statements.
\begin{lemma}\label{converge}
Let $Q(z)  = \sum_{n=2}^{\infty} (n^3-n) a_n z^{n-2} \in
A_{\infty}(\Delta)$, then the series $ \sum_{n=2}^{\infty} n^{2s}
|a_n|^2$ is convergent for all real $s <1$.
\end{lemma}

\begin{proof}
$Q \in A_{\infty}(\Delta) = \{\phi \; \text{holomorphic on }
\Delta: \sup_{z\in \Delta} |\phi(z)(1-|z|^2)^2| < \infty\}$
implies that for any $\al < 1$,
\begin{align*}
\iint\limits_{\Delta} \left\vert Q(z)(1-|z|^2)^2\right\vert^2
\\frac{dxdy}{(1- |z|^2)^{\al}} < \infty,
\end{align*}
where $z=x+iy$. This integral is equal to
\[
\pi \sum_{n=2}^{\infty} (n^3-n)^2 \frac{\Gamma(5-\al)
\Gamma(n-1)}{\Gamma(4+n-\al)} |a_n|^2.
\]
Stirling's formula for the gamma function $\Gamma$ implies that
\begin{align*}
\lim_{n \rightarrow \infty} \frac{\Gamma(n-1) (n^3-n)^2}
{\Gamma(4+n-\al) n^{1+\al}} = 1.
\end{align*}
By comparison test, the series
\[
 \sum_{n=2}^{\infty} n^{1+\al} |a_n|^2
\]
is convergent for all $\al < 1$, which implies the assertion.
\end{proof}

\begin{remark}
We have used the idea of Velling \cite{Ve} in the proof of this
theorem.
\end{remark}

\begin{lemma}[\cite{Zyg}]\label{Zygmund}
If the function $f(z) = \mathfrak{a}_0 + \mathfrak{a}_1 z +
\ldots + \mathfrak{a}_n z^n + \ldots $ is holomorphic on $\Delta$
and continuous on $\Delta \cup S^1$, and the series
$\sum_{n} n |\mathfrak{a}_n|^2$ is convergent, then the series
\begin{align*}
\mathfrak{a}_0 + \mathfrak{a}_1 e^{i\theta} + \ldots +
\mathfrak{a}_n e^{in\theta} +\ldots
\end{align*}
converges uniformly to $f(e^{i\theta})$ on $0 \leq \theta \leq 2 \pi$.
\end{lemma}
 Since $u = \sum_{n=1}^{\infty} a_n z^{n+1}$ is holomorphic on
 $\Delta$ and is continuous on $\C$, Lemma \ref{converge}
 (with $s=\frac{1}{2}$) and Lemma \ref{Zygmund} imply that the series
\[
\sum_{n=1}^{\infty} a_n e^{i(n+1) \theta}
\]
converges uniformly to the continuous function $u\vert_{S^1}
(e^{i \theta})$ on the unit circle $S^1$.

 Similar arguments imply that the series
\[
\sum_{n=0}^{\infty} b_n e^{i(1-n) \theta}
\]
converges uniformly to the continuous function
$v\vert_{S^1}(e^{i \theta})$ on $S^1$.

Taking the derivative with respect to $t$ on the relation
$\gamma^t = (g^t)^{-1} \circ f^t$ and setting $t=0$, we have
\begin{align}\label{3relation}
 \dot{\gamma} = - \dot{g} + \dot{f}.
\end{align}
This shows that the series
\[
\sum_{n=1}^{\infty} a_n e^{i(n+1)\theta} - \sum_{n=0}^{\infty}
b_n e^{i(1-n) \theta}
\]
converges uniformly to the function $\dot{\g}\vert_{S^1}$.
In particular it is the Fourier series of $\dot{\g}\vert_{S^1}$.
Let $\mathfrak{u}(e^{i\theta})\frac{\pa}{\pa \theta}$ be
the corresponding vector field, so that $\dot{\g}= iz\mathfrak{u}(z)$
on $S^1$. We have proved that the Fourier series of
$\mathfrak{u}(e^{i\theta})$
\[
\sum_{n\in \Z} c_n e^{in\theta} , \hspace{1.5cm} c_{-n} = \ov{c_n}
\]
converges uniformly to $\mathfrak{u}(e^{i\theta})$ . Moreover,
\begin{align*}
i \sum_{n\in Z} c_n e^{i(n+1)\theta} = \sum_{n=1}^{\infty}
a_n e^{i(n+1)\theta} - \sum_{n=0}^{\infty} b_n e^{i(1-n)\theta}.
\end{align*}
Comparing coefficients, we have
\begin{align}\label{tangentvector1}
a_n &=  i c_n \;, \hspace{1cm} b_n= - i c_{-n} \hspace{1cm} n\geq 1,
\end{align}
Moreover, we have the relation
\begin{align}\label{tangentvector2}
a_n = \ov{b_n}.
\end{align}
\end{proof}

By imposing extra normalization conditions, we can pass from
the models for $\mathcal{T}(1)$ to the models for $T(1)$.
\begin{remark}
In \cite{Nag}, Nag proved a result similar to Theorem \ref{W}
for $T(1)$ by using explicit formulas for $\dot{\g}$ and $\dot{f}$
from the theory of quasiconformal mappings. Here we use
a slightly different approach.
\end{remark}

\begin{remark}\label{Sob}
Lemmas \ref{converge}, \ref{Zygmund} and  Theorem \ref{W} imply that
the tangent vectors at the origin of $\text{Homeo}_{qs}(S^1)/S^1$ have
Fourier series $\sum_{n} c_n e^{in \theta}$ that converges absolutely
and uniformly, and belong to the Sobolev class $H^s$, for all $s<1$.
Here the Sobolev space $H^{s}(S^1)$ is defined as
\[
H^s(S^1) =\left\{ \mathfrak{u} (e^{i\theta}) = \sum_{n\in \Z}
\mathfrak{a}_n e^{i n \theta} : \sum_{n \in \Z} |n|^{2s}
|\mathfrak{a}_n|^2 < \infty\right\}.
\]
In light of Theorem \ref{W}, we call a tangent vector
$u = \sum_{n=1}^{\infty} a_n z^{n+1} \in T_0 \mathcal{\tilde{D}}$
to be in $H^{s}$ if it is the image of a $H^s$ vector
$\sum_{n} c_n e^{in \theta}$ under the map $D_0 \mathcal{W}$.
\end{remark}
Combining Ahlfors-Beurling extension theorem and
Bers embedding, we get the map
\begin{align*}
\mathcal{B} : \text{Homeo}_{qs}(S^1)/\Mob(S^1) &\rightarrow
\left(L^{\infty}(\Delta^*)_1/ \sim \right) \rightarrow
A_{\infty}(\Delta),\\
\g & \mapsto \hspace{1cm} [\mu] \hspace{1cm}  \mapsto
\mathcal{S} (w^{\mu} \vert_{\Delta}),
\end{align*}
where $\g = w_{\mu}\vert_{S^1}$. Our argument above
gives immediately
\begin{theorem}\label{Bmap}
The derivative of the map $\mathcal{B}$ at the origin
is the linear mapping $D_0 \mathcal{B} : T_0
\left(\text{Homeo}_{qs}(S^1)/\Mob(S^1)\right) \rightarrow
A_{\infty}(\Delta)$ given by
\begin{align*}
\sum_{n \neq -1,0,1} c_n e^{in \theta} & \mapsto
i \sum_{n=2}^{\infty} (n^3 -n) c_n z^{n-2}.
\end{align*}
\end{theorem}
\vspace{0.3cm}

\subsubsection{More on complex structures}
The almost complex structure $J$ at the origins of
$\Diff_+(S_1)/S^1$ and $\Diff_+(S^1)/\Mob(S^1)$ is
defined by the linear map $J : T_0 \rightarrow T_0$ given by
\begin{align}\label{complexstructure}
J\upsilon = i \sum_{n}  sgn(n) c_n e^{in\theta}
\frac{\pa}{\pa \theta}, \hspace{1cm}  \text{where} \;\;
\upsilon = \sum_{n} c_n e^{in\theta}\frac{\pa}{\pa \theta}.
\end{align}
See references in \cite{NV}. (Notice that we differ from the
definition in \cite{NV} by a negative sign). By Remark \ref{Sob},
$J$ extends to almost complex structures on
$\text{Homeo}_{qs}(S^1)/S^1$ and $\text{Homeo}_{qs}(S^1)/\Mob(S^1)$.

 In \cite{NV}, Nag and Verjoysky proved that the almost
 complex structure $J$ on $\Diff_+(S^1)/\Mob(S^1)$ is
 integrable and corresponding complex structure is the
 pull-back of the complex structure on $T(1)$, induced by
 the complex structure of $L^{\infty}(\Delta)_1$.
 Adapting their proof to our convention, we immediately see
 that the complex structure $J$ on $\text{Homeo}_{qs}(S^1)/S^1$
 coincides with the complex structure induced from $\mathcal{T}(1)$.

Using this convention, the holomorphic tangent vectors are of the form
\begin{align*}
w =\frac{\upsilon - i J \upsilon}{2} = \sum_{n>0} c_n e^{in\theta}
\end{align*}
and the antiholomorphic tangent vectors are of the form
\begin{align*}
\bar{w}= \frac{\upsilon + i J \upsilon}{2} = \sum_{n<0} c_n e^{in\theta}.
\end{align*}

\subsection{Metrics}
We are interested in homogenuous Hermitian metrics, i.e.
Hermitian metrics that are invariant under the right group action
on the homogenuous spaces of $\text{Homeo}_{qs}(S^1)$.
In \cite{Ki} and \cite{KY}, Kirillov and Yuriev studied \Ka metrics
on $\Diff_+(S^1)/S^1$. It is known that the homogeneuos \Ka metrics
on $\Diff_+(S^1)/S^1$ must be of the form
\begin{align}\label{metricfamily}
&\parallel \upsilon \parallel_{a,b}^2 = \sum_{n > 0} (an^3 + bn) |c_n|^2,
\end{align}
where $\upsilon = \sum_{n \in \Z} c_n e^{i n \theta} \frac{\pa}
{\pa \theta} \in T_0 \Diff_+(S^1)/S^1$.
The metric $\parallel \cdot \parallel_{0,1}$ is called the Kirillov metric.

On the other hand, since the vector fields $e^{-i \theta}
\frac{\pa}{\pa \theta}, \frac{\pa}{\pa \theta}, e^{i \theta}
\frac{\pa}{\pa \theta}$ generate the $\PSU(1,1)$ action on $S^1$,
\eqref{metricfamily} defines a metric on $\Diff_+(S^1)/\Mob(S^1)$
if and only if $an^3 + bn =0$ for $n = -1, 0 ,1$. This implies that
up to a constant, there is a unique homogeneuos \Ka metric on
$\Diff_+(S^1)/\Mob(S^1)$ given by
\begin{align}\label{WPmetric}
\parallel \upsilon \parallel^2 = \frac{\pi}{2} \sum_{n > 0}
(n^3 - n) |c_n|^2.
\end{align}

Let $\Gamma$ be a Fuchsian group realized as a subgroup of
$\PSU(1,1)$ acting on $\Delta^*$. Let $L^{\infty}(\Delta^*, \Gamma)$
be the space of Beltrami differentials for $\Gamma$, i.e.
\begin{align*}
L^{\infty}(\Delta^*, \Gamma) = \left\{ \mu \in L^{\infty}(\Delta^*) :
\mu \circ \g \frac{\ov{\g'}}{\g'} = \mu, \forall \g \in \Gamma \right\} \;.
\end{align*}
The \Te space of $\Gamma$, $T(\Gamma)$ is the subspace of the
universal Teichm\"uller space
\begin{align*}
T(\Gamma) = L^{\infty}(\Delta^*, \Gamma)_1 / \sim \;,
\end{align*}
where
\[
L^{\infty}(\Delta^*, \Gamma)_1 = L^{\infty}(\Delta^*)_1
\cap L^{\infty}(\Delta^*, \Gamma) \; ,
\]
and $\sim$ is the same equivalence relation we use to define $T(1)$.
The tangent space at the origin of $T(\Gamma)$ is identified with
the space of harmonic Beltrami differentials of $\Gamma$
\[
\Omega^{-1,1}(\Delta^*, \Gamma) = \Omega^{-1,1}(\Delta^*)
\cap L^{\infty}(\Delta^*, \Gamma)\; .
\]

When $\Gamma$ is a cofinite Fuchsian group, i.e. when the
quotient Riemann surface $\Gamma \bk \Delta^*$ has finite
hyperbolic area, there is a canonical Hermitian metric on
$T(\Gamma)$ given by
\begin{align*}
\la \mu, \nu \ra = \iint\limits_{\Gamma \bk \Delta^*} \mu \ov{\nu}
\rho \;, \hspace{1.5cm} \mu, \nu \in \Omega^{-1,1}(\Delta^*, \Gamma)\;,
\end{align*}
where $\rho$ is the area form of the hyperbolic metric on
$\Delta^*$. This metric is called Weil-Petersson metric.
The notation $T(1)$ for the universal \Te space indicates
that it corresponds to the case $\Ga=\{\id\}$. This suggests
to define the Weil-Petersson metric on $T(1)$ by
\begin{align*}
\la \mu, \nu \ra = \iint\limits_{\Delta^*} \mu \ov{\nu}
\rho \;, \hspace{1.5cm} \mu, \nu \in \Omega^{-1,1}(\Delta^*).
\end{align*}
However, this integral does not converge for all $\mu, \nu \in
\Omega^{-1,1}(\Delta^*)$. In particular, it diverges when both
$\mu, \nu$ are Beltrami differentials of a Fuchsian group that
contains infinitely many elements. However, it is proved by Nag
and Verjoysky in \cite{NV} that the integral is convergent on the
Sobolev class $H^{\frac{3}{2}}$ vector fields, which contains the
$C^2$ class vector fields. More precisely, they proved that the
pull back of the Weil-Petersson metric on $T(1)$ to
$\Diff_+(S^1)/\Mob(S^1)$ coincides with the unique homogenuous \Ka
metric \eqref{WPmetric} on $\Diff_+{S^1}/\Mob(S^1)$ (up to a
factor $4$). Henceforth, when we say the Weil-Petersson metric on
$T(1)$, we understand that it is only defined on tangent vectors
in the Sobolev class $H^{\frac{3}{2}}$.

Under the Bers embedding, the Weil-Petersson metric on
$T(1)$ induces a metric on $A_{\infty}(\Delta)$. It is given by

\begin{theorem}\label{WPmetric4}
For $Q = u_{zzz} \in A_{\infty}(\Delta)$, identified as
a tangent vector to $T(1)$ at the origin such that
$u = \sum_{n=1}^{\infty} a_n z^{n+1} \in H^{\frac{3}{2}}$,
the Weil-Petersson metric has the following form
\begin{align*}
\parallel Q \parallel_{WP}^2 = \frac{\pi}{2} \sum_{n=2}^{\infty}
(n^3 - n) |a_n|^2 = \frac{1}{4} \iint\limits_{\Delta}
|Q(z)|^2 (1-|z|^2)^2 dxdy
\end{align*}
\end{theorem}
\begin{proof}
The first equality follows immediately from the identification
of tangent spaces given by Theorem \ref{Bmap}. The second equality
is an explicit computation of the integral.
\end{proof}

\begin{remark}
The derivative of the map $\tilde{\mathcal{D}} \hookrightarrow
A_{\infty}(\Delta)$ at the origin, $\dot{f} \mapsto \dot{f}_{zzz}$
can be viewed as a linear mapping sending vector fields to quadratic
differentials. The theorem states that the Weil-Petersson metric on
$A_{\infty}(\Delta)$ given by the Bers embedding $T(1) \hookrightarrow
A_{\infty}(\Delta)$ is the usual Weil-Petersson metric defined
on the space of quadratic differentials. This can also be proved
directly by using the isomorphism \eqref{Bersmap}. In particular, we have
\begin{align*}
\parallel Q \circ \g (\g')^2 \parallel_{WP}^2 =
\parallel Q \parallel_{WP}^2 , \hspace{1cm} \text{for all}
\;\; \g\in \PSU(1,1).
\end{align*}
\end{remark}

\begin{remark}
Analogs of Theorems \ref{W}, \ref{Bmap} and \ref{WPmetric4} hold
for finite dimensional \Te spaces $T(\Ga)$ embeded in the
universal \Te space $T(1)$.
\end{remark}

According to Remark \ref{Sob}, Kirillov metric on $\Diff_+(S^1)/S^1$
extends to $\mathcal{T}(1)$. Namely, at the origin, it is of the form
\begin{align}\label{Kimetric}
\parallel \upsilon \parallel^2 = \sum_{n > 0} n |c_n|^2 ,
\end{align}
where $\upsilon = \sum_{n} c_n e^{in\theta} \frac{\pa}{\pa \theta}$
is the corresponding tangent vector. The series \eqref{Kimetric} is
convergent. Using the right translations, we define a homogenuous
\Ka metric on $\mathcal{T}(1)$.

 Since every homogenuous \Ka metric on $Diff_+(S^1)/S^1$ can be
 written as a linear combination of the metric \eqref{Kimetric} and the
 Weil-Petersson metric, and only the former is convergent for
 all the tangent vectors of $\mathcal{T}(1)$, we have
\begin{theorem}
Every homogenuous \Ka metric on $\mathcal{T}(1)$ is a multiple of
the metric \eqref{Kimetric}.
\end{theorem}

\section{Velling's Hermitian Form and Velling-Kirillov Metric}
\subsection{Spherical area theorem}
Let the spherical area of a domain $\Omega$ in $\hat{\C}$ be
\[
A_S(\Omega) = \iint\limits_{\Omega} \frac{ 4 dxdy}{(1+ |z|^2)^2}.
\]
It is invariant under rotation, i.e. $A_S(\Omega) =
A_S(e^{i \theta} (\Omega))$.

Following Velling \cite{Ve}, for $Q \in A_\infty(\Delta)$ and
$t$ small, we consider the one parameter family of functions
$f^{tQ} \in \mathcal{D}$ satisfying $\mathcal{S}(f^{tQ}) = tQ$
and the spherical areas of the domains $\Omega_t = f^{tQ}(\Delta)$,
\begin{align*}
A_S(\Omega_t) &= \iint\limits_{\Omega_t} \frac{ 4 dxdy}{(1+ |z|^2)^2}\\
  &= 4 \iint\limits_{\Delta} \frac{ |df^{tQ}|^2}{(1+ |f^{tQ}|^2)^2}.
\end{align*}

Velling`s spherical area theorem is the following.
\begin{theorem}[Velling\cite{Ve}]\label{sphericalarea}
For $Q \in A_\infty(\Delta)$, we have
\begin{align*}\label{velling}
\frac{d}{d t} A_S(f^{tQ}(\Delta)) \vert_{t=0} = 0,\\
\frac{d^2}{dt^2} A_S(f^{tQ}(\Delta)) \vert_{t=0} \geq 0.
\end{align*}
with equality if and only if $Q=0$.
\end{theorem}
This follows from another result, proved by applying the classical
area theorem.
\begin{theorem}[Velling\cite{Ve}]\label{Velling}
Let $f: \Delta \longrightarrow \hat{\C}$ be a univalent function
(perhaps meromorphic) such that it has Taylor expansion
$f(z)= z(1+ a_2 z^2+ a_3 z^3 + \ldots)$ at the origin.
Then the spherical area $A_S (f(\Delta))$ satisfies
\[
A_S(f(\Delta)) \geq 2 \pi,
\]
with equality if and only if $f= \id$.
\end{theorem}

The second inequality in Velling's spherical area theorem
implies that $\frac{d^2}{dt^2} A_S(f^{tQ}(\Delta)) \vert_{t=0}$
is a Hermitian form on $A_{\infty}(\Delta)$.
Our goal is to compute this form explicitly.

The following lemma is very useful for computations.
\begin{lemma}[\cite{Zyg}]\label{com}
Let $f(z)= \sum_{n=0}^{\infty} a_n z^n$ be an analytic
function on $\Delta$ and $\phi(r)$ an integrable function on $[0,1)$.
Then
\begin{align*}
\iint\limits_{\Delta} \phi(|z|)\text{Re}\,(f(z)) dxdy &=
2 \pi \text{Re}\,(a_0) \int_0^1 \phi(r) dr,\\
\iint\limits_{\Delta} \phi(|z|) |f(z)|^2 dxdy &=
2 \pi \sum_{n=0}^{\infty} |a_n|^2 \int_0^1 \phi(r) r^{2n+1} dr.
\end{align*}
\end{lemma}

\subsection{Velling's Hermitian form}

Now we compute Velling's Hermitian form
$\frac{d^2}{dt^2} A_S(f^{tQ}(\Delta)) \vert_{t=0}$.
For $t$ small, we write the perturbative expansions
\begin{align}
f^{tQ}(z) &= z + tu(z) + t^2 v(z) + O(t^3),\\
u(z) &= z( a_2 z^2 + a_3 z^3 + \ldots) =
\sum_{n=2}^{\infty} a_n z^{n+1},\label{series}\\
v(z) &= z( b_2 z^2 + b_3 z^3 + \ldots) =
\sum_{n=2}^{\infty} b_n z^{n+1}.
\end{align}
Taking $t$ derivative of the equation
$\mathcal{S}(f^{tQ}) = tQ$ and setting $t=0$,
we get the following relation
\begin{align*}
\frac{\pa^3}{\pa z^3} u(z) =  Q(z), \hspace{0.8cm} \text{i.e.} \;\;
Q(z) = \sum_{n=2}^{\infty} (n^3-n) a_n z^{n-2}.
\end{align*}

Using the expansion
\begin{align*}
\frac{ |f_z^{tQ}|^2 }{ (1+|f^{tQ}|^2)^2} &=
\frac{ |1+ tu_z + t^2 v_z|^2}{ (1+ |z + tu + t^2 v|^2)^2 } + O(t^3),\\
\end{align*}
we get
\begin{align*}
\frac{d^2}{dt^2} A_S(f^{tQ}(\Delta))\bigr\vert_{t=0} &=
8 \iint\limits_{\Delta} \frac{\chi(z) dxdy}{(1 + |z|^2)^2},  \\
\chi(z) = (v_z + \ov{v_z}  + |u_z|^2) &-
2\frac{ z\ov{v} + \z v + |u|^2 + (z \ov{u} + \z u)
(u_z + \ov{u_z})}{1+|z|^2} + 3 \frac{(z\ov{u}+ \z u)^2}{ (1+ |z|^2)^2}.
\end{align*}
Using the series expansion ~\eqref{series} and $v(0)=v'(0)=0$,
we see that $v$ drops out from the integration. Applying Lemma
\ref{com} we get
\begin{align}\label{Vere}
\frac{d^2}{dt^2} A_S(f^{tQ}(\Delta))\bigr\vert_{t=0} =
16 \pi \sum_{n=2}^{\infty} \mathfrak{c}_n |a_n|^2,
\end{align}
where
\begin{align*}
\mathfrak{c}_n = \int_0^1  \left( \frac{6r^{2n+4}}{(1+r^2)^4}
- \frac{(4n+6) r^{2n+2}}{(1+r^2)^3} +
\frac{(n+1)^2 r^{2n}}{(1+r^2)^2}\right) r dr.
\end{align*}
We compute $\mathfrak{c}_n$ by repeatedly using
integration by parts:
\begin{align*}
\mathfrak{c}_n &= \frac{1}{2} \int_0^1 \left( \frac{6r^{n+2}}{(1+r)^4}
- \frac{(4n+6) r^{n+1}}{(1+r)^3} + \frac{(n+1)^2 r^n}{(1+r)^2}\right)  dr ,\\
\int_0^1 \frac{r^{n+2}}{(1+r)^4} dr &= -\frac{2n^2 + 7n +7}{24}
+ \frac{n(n+1)(n+2)}{6} \int_0^1 \frac{r^{n-1}}{1+r} dr ,\\
\int_0^1 \frac{r^{n+1}}{(1+r)^3} dr &= -\frac{2n+3}{8}
+ \frac{n(n+1)}{2} \int_0^1 \frac{r^{n-1}}{1+r} dr ,\\
\int_0^1 \frac{r^{n}}{(1+r)^2} dr &= -\frac{1}{2} +
n \int_0^1 \frac{r^{n-1}}{1+r} dr.
\end{align*}
Substituting into $\mathfrak{c}_n$, all the terms
with integrals cancel and we are left with
\[
\mathfrak{c}_n = \frac{n}{8}.
\]
Therefore, we have
\begin{theorem}\label{sphericalmetric}
Let $Q \in A_{\infty}(\Delta)$, then
\begin{align*}
\frac{d^2}{dt^2} A_S(f^{tQ}(\Delta))\Bigr\vert_{t=0} =
2 \pi \sum_{n=2}^{\infty} n |a_n|^2.
\end{align*}
\end{theorem}
Remark \ref{Sob} implies that the series is convergent
for all $Q \in A_{\infty}(\Delta)$.
Hence, we can define a Hermitian form on $A_{\infty}(\Delta)$ by
\begin{align*}
\parallel Q \parallel_S^2 &= \frac{1}{2 \pi} \frac{d^2}{dt^2}
A_S(f^{tQ}(\Delta)) \vert_{t=0} = \sum_{n=2}^{\infty} n |a_n|^2, \\
\text{where} \quad Q(z) &= \sum_{n=2}^{\infty} (n^3-n) a_n z^{n-2},
\end{align*}
which we call Velling's Hermitian form.

\begin{remark}
The first half of the computation above is reproduced from
Velling's unpublished manuscript \cite{Ve}. Velling gave the
result in terms of \eqref{Vere}. Our observation is that
$\mathfrak{c}_n$ can be computed explicitly.
\end{remark}

Notice that in evaluating the Hermitian form, we have chosen a
particular normalized solution $f^{tQ}$ to the equation
$\mathcal{S}(f^{tQ})= tQ$. Any other choices will differ from this
one by post-composition with a $\PSL(2, \C)$ transformation.
However the spherical area of a domain $A_S(f(\Delta))$ is not
invariant if $f$ is post-composed with a $\PSL(2, \C)$
transformation. If we choose different normalization conditions to
identify $T(1)$ as a subgroup of $\mathcal{T}(1)$, we get a
different right invariant metric on $T(1)$. Hence the Hermitian
form $\parallel \cdot \parallel_S$ does not naturally define a
right invariant metric on $T(1)$.

On the other hand, since the correspondence between
$\g \in \text{Homeo}_{qs}(S^1)/S^1$ and $f \in \tilde{\mathcal{D}}$
 is canonical, we can use the same approach to define a
 metric on $\mathcal{T}(1) \simeq \text{Homeo}_{qs}(S^1)/S^1$.
 Namely, given the tangent vector $\upsilon =
 \sum_{n \neq 0} c_n e^{i n \theta}\frac{\pa}{\pa \theta}$
 at the origin with the associated one parameter flow $\gamma^t
 = (g^t)^{-1} \circ f^t \vert_{S^1}$, we define a Hermitian form by
\begin{align*}
\parallel \upsilon \parallel^2 = \frac{1}{2\pi}
\frac{d^2}{dt^2}\bigr\vert_{t=0}  A_S(f^t(\Delta)).
\end{align*}
The proof above holds with an extra term $n=1$
(notice that we only need the fact there are no constant
terms and terms linear in $z$ in the first and
second order perturbations), and we get
\[
\parallel \upsilon \parallel^2= \frac{1}{2\pi}
\frac{d^2}{dt^2}\bigr\vert_{t=0}  A_S(f^t(\Delta))
= \sum_{n=1}^{\infty} n |a_n|^2 = \sum_{n=1}^{\infty} n|c_n|^2,
\]
which coincides with the metric \eqref{Kimetric} at the origin. It
is quite remarkable that this metric, introduced by Velling using
classical function theory coincides with the metric introduced by
Kirillov using orbit method. Henceforth, we call this metric on
$\mathcal{T}(1)$ the Velling-Kirillov metric.

\section{Metrics on \Te Spaces}
\subsection{Universal \Te space}

Let $\kappa$ be the symplectic form of the Velling-Kirillov metric
on $\mathcal{T}(1) \simeq \tilde{\mathcal{D}}$. We want to define
a metric on $T(1)$ by vertical integration of the $(2,2)$ form
$\kappa \wedge \kappa$. Namely, let
\[
\omega = \iint\limits_{\text{fiber}} \kappa \wedge \kappa
\]
and define a Hermitian metric on $T(1)$ such that $\omega$ is  the
corresponding symplectic form \footnote{Since the fiber is not
compact, it is not a priori clear that we will get a well-defined
symplectic form on $T(1)$.}. Since $\kappa$ defines a right
invariant metric, $\omega$ also defines a right invariant metric.
Hence we only have to compute the form $\omega$ at the origin of
$T(1)$. We identify the tangent space of $\mathcal{T}(1)$ at the
origin with $A_{\infty}(\Delta)\oplus \C$. The vertical tangent
space is spanned by $\frac{\pa}{\pa w}$ and $\frac{\pa}{\pa \w}$,
where $w$ is the coordinate on $\C$. Observe that the horizontal
and vertical tangent spaces are orthogonal with respect to the
Velling-Kirillov metric. Hence given a holomorphic tangent vector
$Q \in A_{\infty}(\Delta)$, we have
\begin{align*}
\omega( Q, \bar{Q}) &= 2 \iint\limits_{\Delta^*} \kappa(\hat{Q},
\bar{\hat{Q}}) \kappa(\frac{\pa}{\pa w}, \frac{\pa}{\pa \w}) dw \wedge d\w,\\
\end{align*}
where $\hat{Q}$ is the horizontal lift of $(Q,0)$ to every point
on the fiber. Using the right invariance of the Velling-Kirillov
metric, it immediately follows that $\kappa(\frac{\pa}{\pa w},
\frac{\pa}{\pa \w}) dw \wedge d\w$ is the area form of a right
invariant metric on $\Delta^*$. Hence up to a constant, it is the
hyperbolic area form $dA_H$. Checking at the origin, we find that
\[
\kappa(\frac{\pa}{\pa w}, \frac{\pa}{\pa \w}) dw \wedge d\w =
\frac{ dxdy}{(1-|w|^2)^2}= \frac{1}{4}dA_H ,\hspace{0.5cm} w = x+iy.
\]
Using the identification \eqref{vectoridentification} and the
isomorphism \eqref{Bersmap}, $\hat{Q}$ at $(0,w)$ is identified
with $Q \circ \si_w^{-1} ((\si_w^{-1})')^2$ at the origin. Hence
\[
\kappa(\hat{Q}, \bar{\hat{Q}})(w) = \frac{i}{2} \parallel Q \circ
\si_w^{-1} ((\si_w^{-1})')^2 \parallel_S^2.
\]
Making a change of variable $w \mapsto \frac{1}{w}$, then
$\si_w^{-1}$  is changed to $\g_w$, where modulo $S^1$, $\g_w (z)=
\frac{z+w}{1+ z\w}$. Since pre-composing $Q$ with a rotation does
not change the Hermitian form $\parallel Q \parallel_S^2$, we
finally get
\[
\omega(Q, \bar{Q}) = \frac{i}{4} \iint\limits_{\Delta} \parallel
Q_w  \parallel_S^2 dA_H, \hspace{0.5cm} Q_w = Q \circ \g_w
(\g_w')^2.
\]
Thus our approach to define a Hermitian metric on $T(1)$ coincides
with Velling's suggestion \cite{Ve} of averaging the Hermitian
form $\parallel \cdot \parallel_S^2$ along the fiber to define
Hermitian metric on $T(1)$, i.e.
\begin{align}\label{vemetric}
\parallel Q \parallel_V^2 = -2 i \omega(Q, \bar{Q})
= \frac{1}{2}\iint\limits_{\Delta} \parallel Q_w \parallel_S^2
\frac{4 dxdy}{(1-|w|^2)^2}.
\end{align}

\begin{remark}
I am grateful to my advisor L.Takhtajan for his suggestion of
using vertical integration to obtain a metric on $T(1)$.
\end{remark}

Since the Hermitian form $\parallel Q \parallel_S^2$ is expressed
in terms of the norm square of the corresponding coefficients
$|a_n|^2$, to compute \eqref{vemetric}, it is sufficient to
average $|a_n|^2$ for $n\geq 2$.

We set
\[
Q_w(z) = Q\circ \g_w (\g_w')^2(z) = \sum_{n=2}^{\infty} (n^3-n)
a_n^w z^{n-2}, \hspace{0.5cm} \g_w(z) =\frac{z+w}{1+\ov{w}z}.
\]
Then
\begin{align}\label{coefficient}
a_n^w &=  \frac{1}{(n^3-n)} \frac{(Q\circ \g_w
(\g_w')^2)^{(n-2)}}{(n-2)!}(0),
\end{align}
and
\begin{align*}
\parallel Q_w \parallel_S^2 = \sum_{n=2}^{\infty} n |a_n^w|^2.
\end{align*}

\begin{theorem}\label{averagecoefficient}
Let $u(z) = \sum_{n=1}^{\infty} a_n z^{n+1} \in H^{\frac{3}{2}}$
and $Q = u_{zzz}$. Then
\begin{align*}
\iint\limits_{\Delta}|a_j^w|^2 \frac{4 dxdy}{(1-|w|^2)^2}  &=
\frac{2}{3(j^3 - j)} \iint\limits_{\Delta} |Q(w)|^2 (1-|w|^2)^2
dxdy \\ & =\frac{4\pi}{3(j^3 - j)} \sum_{n=2}^{\infty} (n^3-n)
|a_n|^2.
\end{align*}
\end{theorem}

\begin{proof}
Using ~\eqref{coefficient}, we set
\begin{align*}
a_j^w &=  \frac{1}{(j^3-j)}\frac{(Q\circ \g_w (\g_w')^2)^{(j-2)}}{(j-2)!}(0)
=\frac{c_j (w)}{(j^3-j)},
\end{align*}
and introduce the generating function for the $c_j(w)$' s,
\begin{align*}
f(u,w) &= \sum_{j=2}^{\infty} c_j(w) u^{j-2} \\
&= \sum_{j=2}^{\infty} \frac{(Q\circ \g_w (\g_w')^2)^{(j-2)}}
{(j-2)!}(0)u^{j-2}\\
&= Q\circ \g_w(u) (\g_w'(u))^2.
\end{align*}

Writing $u=\rho e^{i\al}$, we have
\begin{align*}
\sum_{j=2}^{\infty} |c_j(w)|^2 \rho^{2j-4} = \frac{1}{2\pi}
\int_0^{2 \pi} |f(\rho e^{i\al},w)|^2 d\al
\end{align*}
and
\begin{align*}
\sum_{j=2}^{\infty} \iint\limits_{\Delta} |c_j(w)|^2
& \frac{ dxdy}{(1-|w|^2)^2} \rho^{2j-4} = \frac{1}{2\pi} \int_0^{2 \pi}
\iint\limits_{\Delta}|f(\rho e^{i\al},w)|^2 \frac{ dxdy}{(1-|w|^2)^2} d\al \\
&= \frac{1}{2\pi} \int_0^{2 \pi} \iint\limits_{\Delta}|
Q\circ \g_w(\rho e^{i\al}) (\g_w'(\rho e^{i\al}))^2|^2
\frac{ dxdy}{(1-|w|^2)^2} d\al. \\
\end{align*}
Denoting this integral by $\mathcal{I}$, substituting the  series
expansion of $Q$ and using polar coordinates $w=re^{i \theta}$, we
get
\begin{align}\label{integral}
&\mathcal{I} = \frac{1}{2\pi} \int_0^{2 \pi} \int_0^1 \int_0^{2\pi}
d\theta rdr d\al \; \Bigl\vert\frac{(1-r^2)}{(1+ r \rho
e^{i(\al-\theta)})^4}\Bigr\vert^2 \\
&\sum_{n=2}^{\infty} (n^3-n) a_n \left(\frac{\rho e^{i\al} +
re^{i\theta}}{1+ r \rho e^{i(\al-\theta)}}\right)^{n-2}
\sum_{m=2}^{\infty} (m^3-m) \ov{a_m} \left(\frac{\rho e^{-i\al} +
re^{-i\theta}}{1+ r \rho e^{-i(\al-\theta)}}\right)^{m-2} .
\end{align}
We do some ``juggling''
\begin{align*}
 &\left(\frac{\rho e^{i\al} + re^{i\theta}}
 {1+ r \rho e^{i(\al-\theta)}}\right)^{n-2} \left(\frac{\rho e^{-i\al} +
 re^{-i\theta}}{1+ r \rho e^{-i(\al-\theta)}}\right)^{m-2} \\
& = \left(\frac{\rho e^{i(\al-\theta)} + r}
{1+ r \rho e^{i(\al-\theta)}}\right)^{n-2}
\left(\frac{\rho e^{-i(\al-\theta)} + r}
{1+ r \rho e^{-i(\al-\theta)}}\right)^{m-2} e^{i(n-m)\theta}
\end{align*}
and make a change of variable $\al \mapsto (\al + \theta)$ to get
\begin{align*}
\mathcal{I} =& \frac{1}{2\pi} \int_0^{2 \pi} \int_0^1 \int_0^{2\pi}
d\theta rdr d\al \; \sum_{n,m \geq 2} (n^3-n)(m^3-m)a_n \ov{a_m} \\
&\left(\frac{\rho e^{i\al} + r}{1+ r \rho e^{i\al}}\right)^{n-2}
\left(\frac{\rho e^{-i\al} + r}{1+ r \rho e^{-i\al}}\right)^{m-2}
\Bigl\vert\frac{(1-r^2)}{(1+ r \rho e^{i\al})^4}\Bigr\vert^2
e^{i(n-m)\theta} \\
=&\int_0^{2 \pi} \int_0^1 rdr d\al \;\\
& \sum_{n=2}^{\infty} (n^3-n)^2 |a_n|^2 \left(\frac{\rho e^{i\al} + r}
{1+ r \rho e^{i\al}}\right)^{n-2}\left(\frac{\rho e^{-i\al} + r}
{1+ r \rho e^{-i\al}}\right)^{n-2} \Bigl\vert\frac{(1-r^2)}
{(1+ r \rho e^{i\al})^4}\Bigr\vert^2 \\
=& \int_0^{2 \pi} \int_0^1 rdr d\al \;\\
&\sum_{n=2}^{\infty} (n^3-n)^2 |a_n|^2 \left(\frac{\rho +re^{i\al}}
{1+ \rho r e^{i\al}}\right)^{n-2}\left(\frac{\rho + re^{-i\al}}
{1+  \rho re^{-i\al}}\right)^{n-2} \Bigl\vert\frac{(1-r^2)}
{(1+  \rho re^{i\al})^4}\Bigr\vert^2  \\
=& \iint\limits_{\Delta} \sum_{n=2}^{\infty} (n^3-n)^2 |a_n|^2
\left(\frac{\rho + w}{1+ \rho w}\right)^{n-2}\left(\frac{\rho +
\ov{w}}{1+  \rho \ov{w}}\right)^{n-2} \Bigl\vert\frac{1-|w|^2}
{(1+  \rho w)^4}\Bigr\vert^2 dxdy,
\end{align*}
where we have done another juggling to get the second
to the last
equality. Observe that
\begin{align*}
\frac{\rho + w}{1+ \rho w} & = \g_{\rho} (w),\\
\frac{1}{(1+  \rho w)^4} &= \frac{\g_{\rho}' (w)^2}{(1-\rho^2)^2}.
\end{align*}
Hence we have
\begin{align*}
& \iint\limits_{\Delta} \left(\frac{\rho + w}{1+ \rho w}
\right)^{n-2}\left(\frac{\rho + \ov{w}}{1+  \rho \ov{w}}
\right)^{n-2} \Bigl\vert\frac{1-|w|^2}{(1+  \rho w)^4}
\Bigr\vert^2 dxdy \\
=& \iint\limits_{\Delta} \left((z^{n-2})\circ \g_{\rho}
(\g_{\rho}')^2 \right)(w)\ov{\left((z^{n-2})\circ \g_{\rho}
(\g_{\rho}')^2 \right)} (w)\frac{(1-|w|^2)^2}{(1-\rho^2)^4} dxdy\\
=& \iint\limits_{\Delta} w^{n-2}\ov{w^{n-2} }\frac{(1-|w|^2)^2}
{(1-\rho^2)^4} dxdy\\
\end{align*}
using $\PSU(1,1)$--invariance of the Weil-Petersson metric. This gives
\begin{align*}
\mathcal{I} =&  \iint\limits_{\Delta} \sum_{n=2}^{\infty}
(n^3-n)^2 |a_n|^2 w^{n-2}\ov{w^{n-2} }\frac{(1-|w|^2)^2}
{(1-\rho^2)^4} dxdy\\
=& \frac{1}{(1-\rho^2)^4} \iint\limits_{\Delta} |Q(w)|^2
(1-|w|^2)^2 dxdy \\
=&\sum_{j=2}^{\infty} \frac{j^3-j}{6} \rho^{2j-4}
\iint\limits_{\Delta} |Q(w)|^2 (1-|w|^2)^2 dxdy.
\end{align*}
Comparing coefficients, we get
\begin{align*}
\iint\limits_{\Delta} |c_j(w)|^2   \frac{ dxdy}{(1-|w|^2)^2}
&=\frac{j^3-j}{6} \iint\limits_{\Delta} |Q(w)|^2 (1-|w|^2)^2 dxdy, \\
\iint\limits_{\Delta}|a_j^w|^2 \frac{4 dxdy}{(1-|w|^2)^2}
&= \frac{2}{3(j^3-j)}\iint\limits_{\Delta} |Q(w)|^2 (1-|w|^2)^2 dxdy,
\end{align*}
which finishes the proof.
\end{proof}

\begin{theorem}\label{Vellingmetric}
Let $Q = u_{zzz} \in A_{\infty}(\Delta)$ be a tangent vector
to $T(1)$ at the origin such that $u \in H^{\frac{3}{2}}$. Then
\begin{align*}
\parallel Q \parallel_{V}^2 = \iint\limits_{\Delta} \parallel
Q_w \parallel_S^2 \frac{2 dxdy}{(1-|w|^2)^2} =
\frac{1}{4} \iint\limits_{\Delta} |Q(w)|^2 (1-|w|^2)^2 dxdy,
\end{align*}
which is the Weil-Petersson metric.
\end{theorem}

\noindent
\begin{proof} This is just a simple sum of the telescoping series:
\begin{align*}
\iint\limits_{\Delta} \parallel Q_w \parallel_S^2 \frac{2 dxdy}
{(1-|w|^2)^2} &= \sum_{j=2}^{\infty} j \iint\limits_{\Delta}
|a_j^w|^2 \frac{2 dxdy}{(1-|w|^2)^2}\\
&=\sum_{j=2}^{\infty} \frac{1}{3(j-1)(j+1)}\iint\limits_{\Delta}
|Q(w)|^2 (1-|w|^2)^2 dxdy \\
&=\frac{1}{4} \iint\limits_{\Delta} |Q(w)|^2 (1-|w|^2)^2 dxdy.
\end{align*}
\end{proof}

\subsection{Finite dimensional \Te spaces}
Let $\Gamma$ be a cofinite Fuchsian group. The tangent space
to $T(\Gamma)$ at the origin is identified with
\begin{align*}
A_{\infty}(\Delta, \Gamma) = \left\{ Q \in A_{\infty}(\Delta)
: Q\circ \g (\g')^2 = Q, \forall \g \in \Gamma\right\}
\end{align*}
and the Weil-Petersson metric is given by
\begin{align}\label{WPmetric2}
\parallel Q \parallel_{WP}^2 = \frac{1}{4}
\iint\limits_{\Gamma\bk \Delta} |Q(w)|^2 (1-|w|^2)^2 dxdy.
\end{align}

The inverse image of $T(\Gamma)$ under the projection map
$\mathcal{T}(1) \rightarrow T(\Gamma)$ is the Bers fiber space
$\mathcal{BF}(\Gamma)$. The quasi-Fuchsian group $\Gamma^{\mu} =
w^{\mu}\circ \Gamma \circ (w^{\mu})^{-1}$ acts on the fiber
$w^{\mu} (\Delta^*)$ at the point $[\mu] \in T(\Gamma)$. The
quotient space of each fiber is a corresponding Riemann surface.
They glue together to form the fiber space $\mathcal{F}(\Gamma)$
over $T(\Gamma)$, which is called the \Te curve of $\Gamma$. First
we have
\begin{lemma}
Let $\Gamma$ be a Fuchsian group. The symplectic form $\kappa$ on
$\mathcal{T}(1)$ restricted to $\mathcal{BF}(\Gamma)$ is
equivariant with respect to the group action on each fiber.
\end{lemma}

\begin{proof}
We only need to check this statement on the fiber at the origin.
The form $\kappa$ restricted to the vertical direction is clearly
equivariant. We are left to verify that if $w \in \Delta^*$,
$\g\in \Gamma$ and $Q \in A_{\infty}(\Delta, \Gamma)$, then
\[
\kappa(\hat{Q}, \ov{\hat{Q}}) (w) = \kappa(\hat{Q}, \ov{\hat{Q}}) (w'),
\]
where $w' = \g(w)$. Notice that the $\PSU(1,1)$ transformation
$\si_{w'} \circ \g \circ \si_w^{-1}$ fixes $\infty$, hence it is
a rotation. Using the fact that the Hermitian form $\parallel Q
\parallel_S^2$ is invariant if $Q$ is pre-composed with a
rotation, we have
\begin{align*}
\parallel Q\circ \si_{w'}^{-1} \left((\si_{w'}^{-1})'\right)^2 \parallel_S^2
&= \parallel \left(Q\circ \g (\g')^2\right)\circ \si_{w}^{-1}
\left((\si_{w}^{-1})'\right)^2 \parallel_S^2 \\&=
\parallel Q\circ \si_{w}^{-1} \left((\si_{w}^{-1})'\right)^2 \parallel_S^2.
\end{align*}

\end{proof}

The lemma implies that $\kappa$ descends to a well-defined
symplectic form on $\mathcal{F}(\Gamma)$. We vertically integrate
the $(2,2)$--form $\kappa \wedge \kappa$ on $\mathcal{F}(\Gamma)$
to define the Hermitian metric on $T(\Gamma)$. Using the same
reasoning as in Section 4.1, we get
\begin{align}
\parallel Q \parallel_V^2 = \frac{1}{2} \iint\limits_{\Gamma \bk \Delta}
\parallel Q_w \parallel_S^2 dA_H, \hspace{1cm} Q \in A_{\infty}(\Delta, \Gamma).
\end{align}
 We want to compute this integral using a regularization
 technique suggested by J. Velling \cite{Ve}.

\begin{theorem}\label{regularization}
Let $\Gamma$ be a cofinite Fuchsian group and  $h \in
L^{\infty}(\Delta)$ be  $\Gamma$--automorphic. We have the
following formula
\begin{align*}
\iint\limits_{\Gamma\bk \Delta} h(w) dA_H = \lim_{r' \rightarrow
1^-}  \frac{ \Area_H(\Gamma\bk \Delta) \iint_{\Delta_{r'}}  h(w)
dA_H }{\iint_{\Delta_{r'}} dA_H},
\end{align*}
where $\Area_H(\Gamma\bk \Delta)$ is the hyperbolic area of the
quotient Riemann surface $\Gamma \bk \Delta$ and $\Delta_{r'} = \{
z : |z| < r' \}$.
\end{theorem}
\begin{proof}
We use the fact that for any $z \in \Delta$, the number of
elements $\gamma \in \Gamma$ such that $\gamma(z)$ is in the disc
$\Delta_{r'}$, is given asymptotically in terms of $r'$ by
\begin{align}\label{lattice}
\frac{1}{\Area_H(\Gamma\bk \Delta)} \iint\limits_{\Delta_{r'}}
dA_H \Bigl(1+o(1)\Bigr), \hspace{1cm} \text{as} \;\; r'
\rightarrow 1^-,
\end{align}
where the $o(1)$ term is uniform for all z in a compact set (see
\cite{Pa}).

Let $F$ be a fundamental domain of $\Gamma$.
Given $E \subset F$, let $E' = \cup_{\gamma \in \Gamma}
\gamma(E)$. Let $\chi_A$ denotes the characteristic function of
the set $A$. Since $\Gamma$ is cofinite, using \eqref{lattice}, we
have
\begin{align*}
\iint\limits_{F} \chi_E dA_H &= \frac{ \Area_H(\Gamma\bk \Delta)
\iint_{\Delta_{r'}} \chi_{E'} dA_H   }{\iint_{\Delta_{r'}} dA_H} +
o(1).
\end{align*}
Here the $o(1)$ term is uniform for all the sets $E \subset F$.
\[
\sup_{w\in \Delta} |h(w)| < \infty ,
\]
standard approximations of $h$ by bounded step functions give our assertion.
\end{proof}

\begin{corollary}
\begin{align*}
\parallel Q \parallel_{WP}^2 = \lim_{r' \rightarrow 1^-}
\frac{ \Area_H(\Gamma\bk \Delta) \iint_{\Delta_{r'}}
|Q(w)|^2 \frac{(1-|w|^2)^2}{4} dxdy }{\iint_{\Delta_{r'}} dA_H}.
\end{align*}
\end{corollary}
\begin{proof}
Take $h(w) = \left\vert Q(w) (1-|w|^2)^2 \right\vert^2$. $Q \in
A_{\infty}(\Delta)$ implies that $h$ is in $L^{\infty}(\Delta)$.
\end{proof}

\begin{lemma}\label{bounded}
Let $Q \in A_{\infty}(\Delta)$. Then
\begin{align*}
\sup_{ w \in \Delta} \parallel Q_w \parallel_S^2 < \infty.
\end{align*}
\end{lemma}
\begin{proof}
Let $\left(Q \circ \g_w (\g_w')^2 \right)(z)= Q_w(z) = \s
um_{n=2}^{\infty} a_n^w z^{n-2} $. The proof of Lemma
\ref{converge} with $\al = 0$, implies that
\begin{align*}
\parallel Q_w \parallel_S^2 = \sum_{n=2}^{\infty} n |a_n^w|^2
< C \iint\limits_{\Delta} \left\vert Q_w (z) (1-|z|^2)^2
\right\vert^2 dxdy, \hspace{1cm} z = x+iy,
\end{align*}
where $C$ is a constant independent of $Q \in A_{\infty}(\Delta)$.
Making a change of variable $z \mapsto \g_w^{-1}(z)$, the integral
on the right hand side is equal to
\begin{align*}
\iint\limits_{\Delta} \left\vert Q(z) (1-|z|^2)^2 \right\vert^2
\left\vert (\g_w^{-1})'(z) \right\vert^2 dxdy
\end{align*}
Since $Q \in A_{\infty}(\Delta)$, the formula (see \cite{Kra})
\begin{align*}
\iint\limits_{\Delta} \left\vert (\g_w^{-1})'(z) \right\vert^2
dxdy = \iint\limits_{\Delta} \frac{(1-|w|^2)^2}{\vert
1-z\bar{w}\vert^4} dxdy = 2 \pi
\end{align*}
concludes the lemma.

\end{proof}

Theorem \ref{regularization} and Lemma \ref{bounded} imply that
our approach to define a Hermitian metric on $T(\Gamma)$ agrees
with J. Velling's original suggestion of using regularized
integrals. Namely, it follows from Theorem \ref{regularization}
and Lemma \ref{bounded},
\begin{corollary}
\begin{align*}
\parallel Q \parallel_{V}^2 = \frac{1}{2} \lim_{r' \rightarrow 1^-}
\frac{ \Area_H(\Gamma\bk \Delta)  \iint_{\Delta_{r'}}
\parallel Q_w \parallel_S^2 dA_H }{\iint_{\Delta_{r'}} dA_H}.
\end{align*}
\end{corollary}

Now we start to compute $\parallel Q \parallel_V^2$. First we have
\begin{theorem}\label{average2}
Let $\Gamma$ be a cofinite Fuchsian group, $Q\in
A_{\infty}(\Delta, \Gamma)$, then
\begin{align*}
\lim_{r' \rightarrow 1^-} \frac{ \Area_H(\Gamma\bk \Delta)
\iint_{\Delta_{r'}} |a_j^w|^2 dA_H }{\iint_{\Delta_{r'}} dA_H} =
\frac{8}{3(j^3-j)} \parallel Q \parallel_{WP}^2.
\end{align*}
\end{theorem}

\begin{proof}
The proof is almost the same as the proof of Theorem
\ref{averagecoefficient}. We have
\begin{align*}
\mathcal{I}&= \sum_{j=2}^{\infty} \iint\limits_{\Delta_{r'}}
|c_j(w)|^2   \frac{ dxdy}{(1-|w|^2)^2} \rho^{2j-4}\\
&=
\iint\limits_{\Delta_{r'}} \sum_{n=2}^{\infty} (n^3-n)^2 |a_n|^2
\left(\frac{\rho + w}{1+ \rho w}\right)^{n-2}\left(\frac{\rho +
\ov{w}}{1+  \rho \ov{w}}\right)^{n-2} \Bigl\vert\frac{1-|w|^2}{(1+
\rho w)^4}\Bigr\vert^2 dxdy.
\end{align*}

Now observe that if $\g \in \PSU(1,1)$ and $Q \in
A_{\infty}(\Delta, \Gamma) $, then $Q\circ \g (\g')^2 \in
A_{\infty}(\Delta, \g^{-1} \Gamma \g)$, and we have
\begin{align*}
(\parallel Q \parallel_{WP}^2 )_{T(\Gamma)} =(\parallel Q\circ \g
(\g')^2 \parallel_{WP}^2 )_{T( \g^{-1}\Gamma \g)}.
\end{align*}
In particular, for any $u= \rho e^{i \al} \in \Delta$, we have
\begin{align*}
\parallel Q \parallel_{WP}^2= \lim_{r' \rightarrow 1^-}
\frac{  \Area_H(\Gamma \bk \Delta) \iint_{\Delta_{r'}}  |(Q\circ
\g_u (\g_u')^2)(w)|^2 \frac{(1-|w|^2)^2}{4} dxdy
}{\iint_{\Delta_{r'}} dA_H},
\end{align*}
since $\Area_H(\Gamma \bk \Delta) = \Area_H(\g_u^{-1}\Gamma \g_u
\bk \Delta)$. It follows that
\begin{align*}
\parallel Q \parallel_{WP}^2= \lim_{r' \rightarrow 1^-}
\frac{  \Area_H(\Gamma \bk \Delta) \frac{1}{2 \pi} \int_0^{2 \pi}
\iint_{\Delta_{r'}}  |(Q\circ \g_u (\g_u')^2)(w)|^2
\frac{(1-|w|^2)^2}{4} dxdy d\al}{\iint_{\Delta_{r'}} dA_H}.
\end{align*}
But we have
\begin{align*}
& \frac{1}{2 \pi} \int_0^{2 \pi} \iint\limits_{\Delta_{r'}}
|(Q\circ \g_u (\g_u')^2)(w)|^2 \frac{(1-|w|^2)^2}{4} dxdy d\al \\
=& \frac{1}{2\pi} \int_0^{2 \pi} \int_0^{r'} \int_0^{2\pi}
d\theta rdr d\al\; \Bigl\vert\frac{(1-r^2)(1-\rho^2)^2}
{2(1+ r \rho e^{i(\theta-\al)})^4}\Bigr\vert^2 \\
&\sum_{n=2}^{\infty} (n^3-n) a_n \left(\frac{\rho e^{i\al}
+ re^{i\theta}}{1+ r \rho e^{i(\theta - \al)}}\right)^{n-2}
\sum_{m=2}^{\infty} (m^3-m) \ov{a_m} \left(\frac{\rho e^{-i\al}
+ re^{-i\theta}}{1+ r \rho e^{-i(\theta - \al)}}\right)^{m-2}.
\end{align*}
This is similar to the integral ~\eqref{integral} with the role of
$\theta$ and $\al$ interchanged, so it equals to
\begin{align*}
&\frac{(1-\rho^2)^4}{4} \iint\limits_{\Delta_{r'}}
\sum_{n=2}^{\infty} (n^3-n)^2 |a_n|^2 \left(\frac{\rho + w}
{1+ \rho w}\right)^{n-2}\left(\frac{\rho + \ov{w}}{1+  \rho \ov{w}}
\right)^{n-2} \Bigl\vert\frac{1-|w|^2}{(1+  \rho w)^4}\Bigr\vert^2 dxdy \\
&=\frac{(1-\rho^2)^4}{4} \mathcal{I}.
\end{align*}
Hence
\begin{align*}
\sum_{j=2}^{\infty} \lim_{r' \rightarrow 1^-}
\frac{ \Area_H(\Gamma \bk \Delta)\iint_{\Delta_{r'}} |c_j(w)|^2
\frac{ dxdy}{(1-|w|^2)^2}} {\iint_{\Delta_{r'}} dA_H} \rho^{2j-4}\\
&= \frac{4}{(1-\rho^2)^4} \parallel Q \parallel_{WP}^2.\\
\end{align*}
Comparing coefficients, we have
\begin{align*}
\lim_{r' \rightarrow 1^-} \frac{ \Area_H(\Gamma \bk \Delta)
\iint_{\Delta_{r'}} |c_j(w)|^2   \frac{ dxdy}{(1-|w|^2)^2}}
{\iint_{\Delta_{r'}} dA_H} = \frac{2(j^3-j)}{3} \parallel Q \parallel_{WP}^2
\end{align*}
and
\begin{align*}
\lim_{r' \rightarrow 1^-} \frac{ \Area_{H}(\Gamma\bk \Delta)
\iint_{\Delta_{r'}} |a_j^w|^2 dA_H }{\iint_{\Delta_{r'}} dA_H}=
\frac{8}{3(j^3-j)} \parallel Q \parallel_{WP}^2.
\end{align*}
\end{proof}

As in the case of Theorem \ref{Vellingmetric}, this immediately implies
\begin{theorem}\label{Velling2}
Let $\Gamma$ be a cofinite Fuchsian group and $Q \in
A_{\infty}(\Delta, \Ga)$ a tangent vector to $T(\Ga)$ at the
origin. Then
\begin{align*}
\parallel Q \parallel_V^2 =  \parallel Q \parallel_{WP}^2.
\end{align*}
\end{theorem}

For a general Fuchsian group $\Gamma$ and $Q \in
A_{\infty}(\Delta, \Gamma)$, we can define
\begin{align*}
\parallel Q \parallel_{V}^2 = \lim_{r' \rightarrow 1^-}
\frac{ \iint_{\Delta_{r'} \cap F(\Gamma)} dA_H \iint_{\Delta_{r'}}
\parallel Q_w \parallel_S^2 dA_H }{\iint_{\Delta_{r'}} dA_H},
\end{align*}
whenever the limit is finite. Here $F(\Gamma)$ is a fundamental
domain of $\Gamma$ on $\Delta$. When $\Gamma$ is the trivial
group, it reduces to integrating over the whole disc, which
coincides with our original definition.

\appendix
\section{Embedding of $\mathcal{T}(1)$}
Consider the Banach space
\[
\mathcal{A}_{\infty}(\Delta) =\left\{\psi \;\text{holomorphic on }
\Delta: \sup_{z\in \Delta} |\psi(z)(1-|z|^2)| < \infty \right\}.
\]
Analogous to the Bers embedding $T(1) \simeq \mathcal{D}
\hookrightarrow A_{\infty}(\Delta)$ (defined in Section 2.1),
which is achieved by the mapping $f \in \mathcal{D} \mapsto
\mathcal{S}(f) \in A_{\infty}(\Delta)$, we prove that there is an
embedding $\mathcal{T}(1) \simeq \tilde{\mathcal{D}}
\hookrightarrow \mathcal{A}_{\infty}(\Delta)$, achieved by the
mapping $f \in \tilde{\mathcal{D}} \mapsto \theta(f)$, where
\[
\theta(f) = \frac{d}{dz}\log f_z = \frac{f_{zz}}{f_z}.
\]
By the classical distortion theorem (see, e.g., \cite{A3}), $f \in
\tilde{\mathcal{D}}$ implies that
\[
\left\vert \frac{f_{zz}}{f_z} - \frac{2 \z}{(1 - |z|^2)} \right
\vert \leq \frac{4}{1-|z|^2}.
\]
Hence $\theta(f) \in \mathcal{A}_{\infty}(\Delta)$, and the map
$\theta:\tilde{\mathcal{D}} \rightarrow
\mathcal{A}_{\infty}(\Delta)$ is well-defined. We claim that this
map is an embedding, and the image contains an open ball.

\begin{lemma}\label{injective}
The map $\theta$ is injective.
\end{lemma}
\begin{proof}
If $f, g \in \tilde{\mathcal{D}}$ are such that $\theta(f) = \theta(g)$, then
\[
\frac{d}{dz} \log f_z = \frac{d}{dz} \log g_z.\] This implies that
$f = c_1 g + c_2$ for some constants $c_1$ and $c_2$. The
normalization conditions $f(0)=g(0)=0$, $f'(0) = g'(0) =1$ (from
the definition of $f,g \in \tilde{\mathcal{D}}$) implies that
$c_1=1, c_2=0$. Hence $f=g$.
\end{proof}

We use the following notations for the sup-norms of
$\mathcal{A}_{\infty}(\Delta)$ and $A_{\infty}(\Delta)$.
\begin{align*}
&\parallel \psi \parallel_{\infty,1} = \sup_{z\in \Delta}
|\psi(z)(1-|z|^2)|, \hspace{0.5cm} \psi \in \mathcal{A}_{\infty}(\Delta);\\
&\parallel \phi \parallel_{\infty,2} = \sup_{z\in \Delta}
|\phi(z)(1-|z|^2)^2|, \hspace{0.5cm}  \phi \in A_{\infty}(\Delta).
\end{align*}
Notice that $\mathcal{S}(f) = \theta(f)_z - \frac{1}{2}
\theta(f)^2$. For $\psi \in \mathcal{A}_{\infty}(\Delta)$, we
define
\[
\psi \mapsto \Psi(\psi) = \psi_z - \frac{1}{2}\psi^2.
\]
We claim that this is a map from $\mathcal{A}_{\infty}(\Delta)$ to
$A_{\infty}(\Delta)$. First, we have the following continuity
theorem.

\begin{theorem}\label{continuity}
For any $\vep >0$, there exists $\delta >0$, such that if $\psi
\in \mathcal{A}_{\infty}(\Delta)$ satisfies $\parallel
\psi\parallel_{\infty,1} < \delta$, then $\Psi(\psi) \in
A_{\infty}(\Delta)$ and $\parallel \Psi(\psi) \parallel_{\infty,2}
< \vep$.
\end{theorem}
\begin{proof}
Fix $\de>0$ and assume that $\psi \in
\mathcal{A}_{\infty}(\Delta)$ satisfies $\parallel \psi
\parallel_{\infty,1} < \delta$. We use the Cauchy formula
\[
\psi_z(z) = \frac{1}{2\pi i} \oint_{|w|=r} \frac{\psi(w)}{(w-z)^2}
dw, \hspace{1cm} |z| < r <1,
\]
to estimate $\psi_z(z)$. Since $\sup_{w \in \Delta} |\psi(w)
(1-|w|^2)| < \delta$,
\[
|\psi_z(z)| \leq \frac{\de}{2\pi (1-r^2)} \oint_{|w|=r} \frac{|dw|}{|w-z|^2}.
\]
Elementary computation gives
\[
\frac{1}{2 \pi} \oint_{|w|=r} \frac{|dw|}{|w-z|^2} = \frac{r}{r^2 - |z|^2}.
\]
Choosing $r = \frac{1+ |z|}{2}$, after some elementary
computations, we have \footnote{This is not the sharpest
estimate.}
\[
|\psi_z(z)(1-|z|^2)^2| \leq \frac{ 8 \de
(1+|z|)^3}{(|z|+3)(1+3|z|)} \leq \frac{64 \de}{3}  \hspace{0.5cm}
\text{for} \;\; |z|\leq 1.
\]
Hence
\begin{align*}
\left\vert \left(\psi_z(z) -
\frac{1}{2}(\psi(z))^2\right)(1-|z|^2)^2\right\vert \leq
\frac{64\de}{3} + \frac{\de^2}{2}.
\end{align*}
Given $\vep >0$, we can always find $\de >0$, such that
$\frac{64\de}{3} + \frac{\de^2}{2} < \vep$. This proves our
assertion.

\end{proof}

\begin{corollary}\label{holomorphy}
$\psi \mapsto \Psi(\psi)$ is a holomorphic map from
$\mathcal{A}_{\infty}(\Delta)$ to $A_{\infty}(\Delta)$.
\end{corollary}
\begin{proof}
The map $\psi \mapsto \psi_z$ is linear. From the proof of the
theorem above, we see that it is a continuous map from
$\mathcal{A}_{\infty}(\Delta)$ to $A_{\infty}(\Delta)$. The map
$\psi \mapsto -\frac{1}{2} \psi^2$ is clearly a continuous map
from $\mathcal{A}_{\infty}(\Delta)$ to $A_{\infty}(\Delta)$. Hence
$\Psi$ is a continuous map from $\mathcal{A}_{\infty}(\Delta)$ to
$A_{\infty}(\Delta)$.

To prove holomorphy, it is sufficient to note that for any $\psi,
\varphi \in \mathcal{A}_{\infty}(\Delta)$ and $\ep \in \C$ in a
neighbourhood of $0$, the Frechet derivative
\begin{align*}
\lim_{\ep \rightarrow 0} \frac{ \Psi(\psi + \ep \varphi) -
\Psi(\psi)}{\ep} = \varphi_z - \psi \varphi
\end{align*}
exists in the $\parallel \cdot \parallel_{\infty, 2}$ norm, since
\begin{align*}
\Bigl\Vert  \frac{ \Psi(\psi + \ep \varphi) - \Psi(\psi)}{\ep} -
(\varphi_z - \psi \varphi) \Bigr\Vert_{\infty,2} =  \ep \parallel
\frac{1}{2} \varphi^2 \parallel_{\infty,2} = \frac{1}{2} \ep
\parallel \varphi \parallel_{\infty,1}^2
\end{align*}
tends to $0$ as $\ep \rightarrow 0$.
\end{proof}

\begin{theorem}\label{openball}
The image of $\theta$ contains an open ball about the origin of
$\mathcal{A}_{\infty}(\Delta)$.
\end{theorem}
\begin{proof}
From Theorem \ref{continuity}, there exists an $\al$ such that if
$\parallel \psi \parallel_{\infty,1} < \al$, then $\phi =
\Psi(\psi)$ satisfies $\parallel \phi \parallel_{\infty,2} < 2$.
By Ahlfors-Weill theorem, there exists a univalent function $f_1:
\Delta \rightarrow \C$ such that $\mathcal{S}(f_1) = \phi$ and
$f_1$ has quasiconformal extension to $\C$. On the other hand,
there exists a unique holomorphic function $f: \Delta \rightarrow
\C$ which solves the ordinary differential equation
\[
\frac{d}{dz} \log  f_z = \psi; \hspace{0.5cm}  f(0)=0, f'(0)=1.
\]

Obviously, $\mathcal{S}(f) = \Psi(\psi) = \phi$. Hence $f$ and
$f_1$ agrees up to post-composition with a $\PSL(2, \C)$
transformation. This implies that $f$ also has a quasiconformal
extension to $\C$ and $f \in \tilde{\mathcal{D}}$. Hence the image
of $\theta$ contains the open ball of radius $\al$.

\end{proof}

From Lemma \ref{injective} and Theorem \ref{openball}, it follows
that there is an embedding of $\tilde{\mathcal{D}}$ into
$\mathcal{A}_{\infty}(\Delta)$ whose image contains an open ball
about the origin. This implies that $\tilde{\mathcal{D}}$ has a
Banach manifold structure modeled on
$\mathcal{A}_{\infty}(\Delta)$. We want to compare this structure
to the structure induced from the embedding $\tilde{\mathcal{D}}
\simeq \mathcal{T}(1) \hookrightarrow A_{\infty}(\Delta)\oplus
\C$. We define a map $\widehat{\Psi} :
\mathcal{A}_{\infty}(\Delta) \rightarrow A_{\infty}(\Delta)\oplus
\C$ by
\[
\psi \mapsto \left(\Psi(\psi), \frac{1}{2} \psi(0)\right).
\]

\begin{theorem}
The map $\widehat{\Psi}$ is holomorphic and one to one.
\end{theorem}

\begin{proof}
Holomorphy follows directly from corollary \ref{holomorphy}. To
prove injectivity, suppose $\widehat{\Psi}(\psi_1) =
\widehat{\Psi}(\psi_2)$. For $j =1,2$, let
\[
f_j (z) = \int_0^{z} e^{\int_0^w \psi_j (u) du} dw.
\]
Then $\frac{d}{dz}\log f_j' = \psi_j$, $f_j(0)=0$, $f_j'(0) =1$.
This implies that $\mathcal{S}(f_1) = \Psi(\psi_1) = \Psi(\psi_2)
= \mathcal{S}(f_2)$. Hence $f_1 = \si \circ f_2$ for some $\si \in
\PSL(2, \C)$. $f_j (0) =0, f_j'(0)=1$, $j=1,2$ implies that $\si =
\ma{1}{0}{c}{1}$ for some $c \in \C$. We also have
\[
\frac{d}{dz} \log f_1' = \frac{d}{dz} \left(\log \si' \circ f_2 +
\log f_2'\right).\] Setting $z=0$, we have
\[
\psi_1(0) = -2c + \psi_2(0).
\]
$\widehat{\Psi}(\psi_1) = \widehat{\Psi}(\psi_2)$ implies $c=0$
and $f_1 = f_2$, $\psi_1 = \psi_2$.

\end{proof}

\begin{theorem}
The Banach spaces $\mathcal{A}_{\infty}(\Delta)$ and
$A_{\infty}(\Delta)\oplus \C$ induce the same Banach manifold
structure on $\tilde{\mathcal{D}}$.
\end{theorem}

\begin{proof}
From our discussion in Section 2.1.2 (in particular
\eqref{isomor}, \eqref{constantc}, \eqref{lambda}), we know that
the embedding $\tilde{\mathcal{D}} \hookrightarrow
A_{\infty}(\Delta)\oplus \C$ factors through the map
$\widehat{\Psi}$, i.e., it is given by $f \mapsto \theta(f)
\xrightarrow{\widehat{\Psi}} (\mathcal{S}(f), \frac{1}{2}
\theta(f)(0))$. Let $U$ (resp. $V$) be the image of
$\tilde{\mathcal{D}}$ in $\mathcal{A}_{\infty}(\Delta)$ (resp.
$A_{\infty}(\Delta)\oplus \C$). Bers proved that $V$ is open in
$A_{\infty}(\Delta) \oplus \C$ (\cite{Bersfiber}) using a theorem
of Ahlfors (\cite{A63}) which says that the image of $T(1)$ in
$A_{\infty}(\Delta)$ is open. The continuity of the map
$\widehat{\Psi}$ implies that $U$ is open in
$\mathcal{A}_{\infty}(\Delta)$. Hence we have a holomorphic
bijection $\widehat{\Psi}\vert_U : U \rightarrow V$. In order to
conclude that this is a biholomorphic map between open subsets of
Banach manifolds, by inverse mapping theorem (see, e.g.,
\cite{Lang}), we only have to show that for any $\psi \in U$, the
derivative of $\widehat{\Psi}$ at $\psi$, $D_{\psi}
\widehat{\Psi}$ is a topological linear isomorphism between
$\mathcal{A}_{\infty}(\Delta)$ and $A_{\infty}(\Delta) \oplus \C$.

From the proof of Corollary \ref{holomorphy}, the linear map
$D_{\psi} \widehat{\Psi} : \mathcal{A}_{\infty}(\Delta)
\rightarrow A_{\infty}(\Delta) \oplus \C$ is given by
\[
D_{\psi}\widehat{\Psi} (\varphi) = \left(\varphi_z - \psi \varphi,
\frac{1}{2} \varphi(0)\right).
\]
From the theory of ordinary differential equations, it is easy to
prove that this map is injective. To prove surjectivity, let $f
\in \tilde{\mathcal{D}}$ be such that $\theta(f) = \psi$. Given
$(\phi, c) \in A_{\infty}(\Delta) \oplus \C$, consider
\[
\varphi(z) = f'(z) \left(\int_{0}^z \frac{\phi(u)}{f'(u)} du + 2c \right).
\]
It is straight-forward to check that $\varphi$ is the unique
holomorphic function on $\Delta$ that satisfies
\[
\varphi_z - \psi \varphi = \phi \hspace{0.5cm} \text{and}
\hspace{0.5cm} \frac{1}{2} \varphi(0) = c.
\]
What remains to be proved is that $\varphi \in \mathcal{A}_{\infty}(\Delta)$.

Let $\Omega = f(\Delta)$ and $\Omega^*$ be the exterior of the
domain $\Omega$. Let $\lambda(w) \vert dw \vert$ be the $\Poin$
metric on $\Omega$, which is given by
\[
\lambda \circ f (z) \vert f'(z)\vert = \frac{1}{(1-|z|^2)}.
\]
For $w \in \Omega$, let $\de(w)$ denotes the Euclidean distance
from $w$ to the boundary of $\Omega$. The Koebe one-quarter
theorem (see, e.g., \cite{Nag2, Lehto}) implies that
\begin{align}\label{onequarter}
\frac{1}{4} \leq \lambda(w) \de(w) \leq 1.
\end{align}
Let $\tilde{\phi}(w) = \phi \circ f^{-1}(w)
\left(f_w^{-1}(w)\right)^2$. Then
\[
\int_0^z \frac{\phi(u)}{f'(u)} du = \int_0^w \tilde{\phi}(v) dv = \Phi(w),
\]
where $w = f(z)$. Since
\[
\sup_{w \in \Omega} \vert \lambda^{-2}(w) \tilde{\phi}(w) \vert
=\sup_{z \in \Delta} \vert (1-|z|^2)^2 \phi(z) \vert  < \infty,
\]
by a theorem of Bers (\cite{Bers66}), there is a bounded harmonic
Beltrami differential on $\Omega^*$, $\mu : \Omega^* \rightarrow
\C$, $\sup_{w \in \Omega^*} \vert \mu(w) \vert = \beta< \infty$
such that
\[
\tilde{\phi}(w) = - \frac{6}{\pi} \iint\limits_{\Omega^*}
\frac{\mu(v)}{(v-w)^4} \left\vert \frac{dv \wedge d\bar{v}}{2}
\right\vert, \hspace{0.5cm} w \in \Omega.
\]
This implies that
\[
\Phi(w) = - \frac{2}{\pi} \iint\limits_{\Omega^*}
\frac{\mu(v)}{(v-w)^3} \left\vert \frac{dv \wedge d\bar{v}}{2}
\right\vert + C_1,
\]
where $C_1$ is a constant so that $\Phi(0) =0$. Since every point
$v \in \Omega^*$ is of distance at least $\de(w)$ away from $w$,
we have the following estimate:
\begin{align*}
\left\vert \Phi(w) \right\vert &\leq \frac{2}{\pi}\iint
\limits_{|v-w| \geq \de(w)} \frac{|\mu(v)|}{|v-w|^3}
\left\vert\frac{dv \wedge d\bar{v}}{2}\right\vert+ C_1\\
& \leq \frac{2 \beta}{\pi} \int_0^{2 \pi}
\int_{\de(w)}^{\infty} \frac{1}{\rho^3} \;\rho d \rho d \theta + C_1\\
&=\frac{4 \beta}{\de(w)} + C_1 \leq 16 \beta \lambda(w) + C_1.
\end{align*}
Using \eqref{onequarter} and $w = f(z)$, we also have
\begin{align*}
\frac{1}{4}  \;\; \leq & \;\;\frac{1}{(1-|z|^2)
\vert f'(z) \vert} \de(f(z))\;\;  \leq \;\;1\\
\hspace{1cm} \frac{\de(f(z))}{1-|z|^2} \;\; \leq &
\hspace{1.5cm}\left\vert f'(z) \right\vert \hspace{1.5cm}  \leq
\;\; \frac{4 \de(f(z))}{1-|z|^2}.
\end{align*}
Hence we have
\begin{align*}
\left\vert \varphi (z)\right\vert &\leq \left\vert f'(z)
\right\vert \left\vert \Phi\left(f(z)\right) + C_2 \right\vert\\
& \leq 16 \beta \left\vert \lambda\left(f(z)\right)\right\vert
\left\vert f'(z)\right\vert +C_2 \frac{4 \de(f(z))}{1-|z|^2}  \\
&\leq \frac{C}{1-|z|^2}.
\end{align*}
Here $C_2$ and $C$ are constants. To get the last inequality, we
have used the fact that $\de(f(z))$ is bounded for $z \in \Delta$.
This concludes that $\varphi \in \mathcal{A}_{\infty}(\Delta)$.

\end{proof}

\begin{remark}
When $\Gamma$ is a Fuchsian group, the embedding $\theta :
\mathcal{\tilde{D}} \simeq \mathcal{T}(1) \rightarrow
\mathcal{A}_{\infty}(\Delta)$ restricts to an embedding
$\mathcal{BF}(\Ga) \rightarrow \mathcal{A}_{\infty}(\Delta, \Ga)$,
where
\[
\mathcal{A}_{\infty}(\Delta, \Ga) = \left\{ \psi \in
\mathcal{A}_{\infty}(\Delta) : \psi_z - \frac{1}{2} \psi^2 \in
A_{\infty}(\Delta, \Ga) \right\}
\]
Contrary to the description of $A_{\infty}(\Delta, \Ga)$ as
$\parallel \cdot \parallel_{\infty, 2}$ bounded holomorphic
quadratic differentials of the Riemann surface $\Gamma \bk
\Delta$, $\mathcal{A}_{\infty}(\Delta, \Ga)$ does not have an
intrinsic characterization as a space of differentials on $\Ga \bk
\Delta$. Rather it is extrinsically defined as the space of
solutions to the Ricatti equation
\[
\psi_z - \frac{1}{2} \psi^2 = \phi, \hspace{1cm} \phi \in A_{\infty}(\Delta, \Ga),
\]
on the Riemann surface $\Ga \bk \Delta$. However, it contains the
subspace of affine connections on $\Ga \bk \Delta$, i.e.
\[
\left\{ \lambda : \Delta \rightarrow \C \;\; \text{holomorphic}\;
: \;\parallel \lambda \parallel_{\infty, 1} < \infty, \lambda\circ
\g \g' - \lambda = \frac{\g''}{\g'} \;\; \forall \g \in \Ga
\right\},
\]
which is an affine space modeled on the space of $\parallel \cdot
\parallel_{\infty, 1}$ bounded holomorphic $1$-forms on $\Gamma\bk
\Delta$.

\end{remark}

\bibliographystyle{amsalpha}
\bibliography{kirillov}

\end{document}